\documentclass[11pt,a4paper,twoside]{article}
\usepackage{bm, amsmath, amssymb, amsthm} 
\usepackage{graphicx}

\newcommand{\ds}{\operatorname{ds}}
\def\<{\langle}
\def\>{\rangle}
\def\eps{\varepsilon}
\def\RR{\mathbb{R}}
\def\calf{\mathcal{F}}
\def\tr{\operatorname{Tr\,}}
\def\id{\operatorname{id\,}}
\def\Div{\operatorname{div}}

\newcommand{\dtau}{\operatorname{d}\!\tau}

\def\Ric{\operatorname{Ric}}
\def\eq{\hspace*{-1.5mm}&=&\hspace*{-1.5mm}}
\def\plus{\hspace*{-1.5mm}&+&\hspace*{-1.5mm}}
\def\minus{\hspace*{-1.5mm}&-&\hspace*{-1.5mm}}
\def\dt{\partial_t}

\newtheorem{corollary}{Corollary}
\newtheorem{definition}{Definition}
\newtheorem{example}{Example}
\newtheorem{remark}{Remark}
\newtheorem{lemma}{Lemma}
\newtheorem{question}{Question}
\newtheorem{conjecture}{Conjecture}
\newtheorem{proposition}{Proposition}
\newtheorem{theorem}{Theorem}

\author{Vladimir Rovenski\thanks{
        Mathematical Department, University of Haifa,
        E-mail: rovenski@math.haifa.ac.il}
        \footnote{The results of Section~\ref{sec:rot_sym} were obtained by V.\,Rovenski in common with L.\,Zelenko.}
   \ and \ Vladimir Sharafutdinov\thanks{
        Sobolev Institute of Mathematics at Novosibirsk,
        E-mail: sharaf@math.nsc.ru}
}

\title{The partial Ricci flow on one-dimensional foliations}

\begin{document}

\date{}

\maketitle

\tableofcontents

\begin{abstract}
A~flow of metrics, $g_t$, on a manifold is a solution of a differential equation $\dt g = S(g)$\,,
where a~geometric functional $S(g)$ is a~symmetric $(0,2)$-tensor usually related to some kind of curvature.
The~mixed sectional curvature of a foliated manifold regulates the deviation of leaves along the leaf geodesics.
We~introduce and study the~flow of metrics on a foliation (called the \textit{Partial Ricci Flow}),
where $S=-2\,r$ and $r$ is the partial Ricci curvature of the foliation; in other words, the velocity
for a unit vector $X$ orthogonal to the leaf, $-2\,r(X,X)$, is the mean value of sectional curvatures over all mixed planes containing~$X$.
 The flow preserves totally geodesic foliations and is used to examine the question:
Which foliations admit a metric with a given property of mixed sectional curvature (e.g., point-wise constant)\/?
This is related to Toponogov question about dimension of totally geodesic foliations with positive mixed sectional curvature.

We first consider a one-dimensional foliation, since this case is easier.
We prove local existence/uniqueness theorem,
deduce the government equations for the curvature and conullity tensors (which are parabolic along the leaves), and show convergence of solution metrics for some classes of almost-product structures.
For the warped product initial metric the global solution metrics converge to one with constant mixed sectional curvature.


\end{abstract}

\section{Introduction}

We discuss Toponogov's question about dimension of totally geodesic foliations with positive mixed sectional curvature (i.e., Ferus's estimate) and define the Partial Ricci Flow on foliations. It exists locally and
is proposed as the main tool to study the question (see the conjecture by~Rovenski).

\subsection{Totally geodesic foliations}

 A Riemannian manifold $(M,g)$ may admit many kinds of geometrically interesting foliations.
The~problems of the existence and classification of metrics on foliations (first posed by H.\,Gluck in 1979 for geodesic foliations) were studied already in the 1970's when D.\,Sullivan provided a~topological condition (called {topological tautness}) for a foliation, equivalent to the existence of a Riemannian metric making all the leaves minimal, see~\cite{cc1}.
Seve\-ral authors investigated whether on a~given Riemannian manifold $(M,g)$ there exists a totally geodesic foliation, as well as the inverse problem of determining whether one can find a Riemannian metric on a foliated manifold $(M,\calf)$ with respect to which the foliation becomes totally geodesic, see for example \cite{G83}, \cite{mk2004} and a survey  in \cite{rov-m}.
 In 1970, Ferus \cite{Fe} proved the following.

\begin{theorem} \label{FT} Let a Riemannian manifold $(M^{n+p}, g)$ be foliated with complete totally geodesic
leaves of dimension $p$. Assume that the sectional curvature of $M$ has the same positive value for all
planes spanned by two vectors such that the first (second) vector is tangent (orthogonal) to a leave.
Then
\begin{equation} \label{FI}
 p\le \rho(n)-1,
\end{equation}
where $\rho(n)-1$ is the maximal number of point-wise linear independent vector fields on a sphere~$S^{n-1}$.
\end{theorem}

To the best of our knowledge, this is the unique theorem in Riemannian geometry, which involves the topological invariant $\rho(n)$, the \textit{Adams number}; here $\rho\big(({\rm odd})\,2^{\,4d+c}\big)=8d+2^{\,c}$\ where $d\ge0$ and $0\le c\le3$.
In the simplest case of $p=1$, the manifold $M$ is foliated to one-dimensional leaves that are complete geodesics. Let us refer to the latter case as a {\it geodesic foliation}.
Let also \textit{mixed curvatures} (\textit{mixed planes}) stand for the sectional curvatures (planes) mentioned in Theorem~\ref{FT}.

Theorem~\ref{FT} prohibits the existence of a geodesic foliation with positive constant mixed curvatures in the case of an even-dimensional manifold, since $\rho(n)-1=0$ for an odd $n$.
Hopf's fiber bundle $\pi:S^3\rightarrow S^2$ gives the simple example of such a foliation for the
odd $n+p=3$, where the sphere $S^3$ is equipped with the standard metric.
Fibers of Hopf's bundle are closed geodesics (great circles).
The theorem has various applications to geometry of submanifolds, for example:
A~complete submanifold $M^l$ of a~sphere $S^{l+m}$ is totally geodesic if
the relative nullity index satisfies the inequality
$\nu(M)\ge F(l) := \max\{s:  s < \rho(l -s)\}$ (the \textit{Ferus number}).
One~may try to find examples and classify submanifolds
in $S^{l+m}$ which satisfy the equality $\nu(M) = F(l)$ (all of them have singularities when considered in $\mathbb{R}^{l+m}$, see recent examples in \cite{ist}).

Among Toponogov's many important contributions to global Riemannian geometry is the following question,
see survey in \cite[p. 30]{rov-m}:

\begin{question}\label{TP}
Can Theorem \ref{FT} be generalized by replacing the hypothesis ``all mixed curvatures are equal to a positive constant'' with the weaker one: ``all mixed curvatures are positive''?
\end{question}

Although the question was posed in 80's, it is still open.
Rovenski \cite{rov-m} proved the exactness of estimate (\ref{FI})
and necessity of more conditions when a foliation is given locally.
He~solved Problem~\ref{TP} for the special case,
when $M^{n+p}$ is a ruled submanifold of a sphere.

\subsection{The partial Ricci flow}

Let $(M^{n+p}, g)$ be a connected Riemannian manifold with the Levi-Civita connection $\nabla$,
$\calf$ a~smooth $p$-dimensional foliation on $M$ and $\mathcal{D}$ its orthogonal $n$-dimensional distribution
(i.e., $g(N,N)=1$ and $g(N,X)=0$ for $X\in\mathcal{D}$ and $N\in\mathcal{D}_\calf$).
 As usual, $R(X,Y,Z,V)=g(R(X,Y)Z,V)$ is the Riemannian curvature tensor, and
$R(X,Y)=\nabla_Y\!\nabla_X{-}\nabla_X\!\nabla_Y{+}\nabla_{[X,Y]}$ is the~curvature operator.
Thus, $R(X,Y)=\nabla^2_{Y,X}-\nabla^2_{X,Y}$ where
$\nabla^2_{X,Y}:=\nabla_{X}{\circ}\nabla_{Y}-\nabla_{\nabla_X Y}$
is the second covariant derivative.

A~flow of metrics, $g_t$, on a manifold is a solution of a differential equation $\dt g = S(g)$\,,
where the symmetric $(0,2)$-tensor $S(g)$ is usually related to some kind of curvature.
The~mixed sectional curvature of a foliated manifold regulates the deviation of leaves along the leaf geodesics.

The first author and Walczak \cite{rw-m} (see also \cite{rw4}) studied flows of metrics
that depend on the extrinsic geometry of codimensi\-on-one foliations.
Rovenski and Wolak~\cite{rovwol} studied ${\mathcal D}$-conformal flows of metrics
on any foliation in order to prescribe the mean curvature vector $H$ of ${\mathcal D}$.

 Denote by $\mathcal{M}$ the space of smooth Riemannian metrics
 on $M$ such that the distribution $\mathcal{D}$ is orthogonal to $\calf$.
 Elements of $\mathcal{M}$ are called $(\mathcal{D},\mathcal{D}_\calf)$-\textit{adapted metrics}
 (\textit{adapted metrics}, in short).

\begin{definition}\rm
The~\textit{Partial Ricci Flow} (PRF) is a family of adapted metrics $g_t,\ t\in[0,\eps)$, satisfying the~PDE
\begin{equation}\label{E-GF-Rmix}
 \dt g = -2\,r(g).
\end{equation}
The symmetric $(0,2)$-tensor $r=r(g)$ (called the \textit{partial Ricci curvature}, \cite{r2010}) is defined as follows:
\begin{equation}\label{E-rF}
 r(X,Y) = \tr^{\calf} R(X^\perp,\,\cdot\,,Y^\perp,\,\cdot\,)\qquad (X,Y\in TM),
\end{equation}
where $^\bot$ is the orthogonal to $\calf$ component of a vector.
\end{definition}

Let $(e_1, \dots, e_p)$ be a local orthonormal basis of the tangent space ${\calf}_x$ to the leaf through a point $x\in M$.
Then $r_x(X,Y)=\sum_{i=1}^p R(X^\bot, e_i, Y^\bot, e_i)$ for $X,Y\in T_xM$.
In~other words, for a unit vector $X$ orthogonal to the leaf, $r(X,X)$ is the mean value of sectional curvatures over all mixed planes containing~$X$.
Observe that $r(X,Y)=0$ if either $X$ or $Y$ is tangent to $\calf$. This means that the PRF does not change geometry of leaves and remains leaves to be totally geodesic.

 The notion of the $\calf$-\emph{truncated} $(0,2)$-\emph{tensor} will be helpful:
\[
 S(X,Y) = S(X^\perp, Y^\perp)\qquad (X,Y\in TM).
\]
The tensor $r$ provides an example of an $\calf$-truncated symmetric $(0,2)$-tensor.
Another useful example is the $\calf$-\textit{truncated metric tensor} $\hat g$,
i.e., $\hat g(X,Y) = g(X^\perp, Y^\perp)$.

\begin{remark}\rm
The trace of $r$ is the \textit{mixed scalar curvature} of a foliation, see \cite{rov-m},
${\rm Sc}_{\,\rm mix}(g)=\tr_g r$.
If $\{e_i,\,\eps_\alpha\}_{i\le p,\,\alpha\le n}$ is a local orthonormal frame on $TM$
adapted to ${\mathcal D}_{\calf}$ and ${\mathcal D}$ then
\[
 {\rm Sc}_{\,\rm mix}(g) =\sum\nolimits_{i=1}^p\sum\nolimits_{\alpha=1}^n R(\eps_\alpha, e_i, \eps_\alpha, e_i).
\]
For a codimension-one foliation with a unit normal $N$, we have ${\rm Sc}_{\,\rm mix}=\Ric(N,N)$.
For a surface $(M^2,g)$, i.e., $n=p=1$, we obtain ${\rm Sc}_{\,\rm mix}=K$ -- the gaussian curvature.
Rovenski and Zelenko \cite{rz2012, rz2013} initiated the study of the mixed scalar curvature flow
\begin{equation}\label{E-GF-Ricmix}
 \dt g = -2\,({\rm Sc}_{\,\rm mix}(g) -\Phi)\,\hat g,
\end{equation}
where $\Phi:M\to{\mathbb R}$ is a leaf-wise constant.
This is the `Yamabe type' analogue to the \textit{normalized}~PRF
\begin{equation}\label{E-GF-Rmix-Phi}
 \dt g = -2\,r(g) +2\,\Phi\,\hat g.
\end{equation}
Metrics with $r=pK\hat g$ for $K\in\RR$ (hence, ${\rm Sc}_{\,\rm mix}=p\,nK$) are the fixed points of \eqref{E-GF-Rmix-Phi}
with $\Phi=pK$.
\end{remark}

Similarly to results of Section~\ref{sec:onedim} one may prove  the following.

\begin{proposition}
The flow \eqref{E-GF-Rmix} preserves the metric of $\mathcal{D}_\calf$ and the orthogonality of vectors to $\calf$. If~$\calf$ is totally geodesic for $t=0$ then it is totally geodesic for all $t$.
\end{proposition}

\begin{theorem}\label{T-main-locF}
Let $(M^{n+p}, g_0)$ be a closed Riemannian manifold with a smooth $p$-dimensional foliation.
Suppose that the leaves compose a fiber bundle with the total space $M$.
Then (\ref{E-GF-Rmix}) has a unique smooth solution $g_t$ defined on a positive time interval~$[0,\eps)$.
\end{theorem}

 The flow may be used to examine the question:
Which foliations admit a metric with a given property of mixed sectional curvature (e.g., point-wise constant)\/?
One may try to attack Question~\ref{TP} by deforming the metric in directions orthogonal to leaves.
The best candidate for such a deformation is the PRF.
We are going to study the PRF along with the same line as the classical Ricci flow is applied in the proof of the smooth $1/4$-pinching sphere theorem, see for example \cite{ah}.
In order to use the Ricci flow machinery, we have to impose some additional restrictions. Indeed, Ricci flow theory is now well-developed for compact manifolds only. There are many open questions in the case of the Ricci flow on a non-compact manifold, most of them are related to the maximum principle for parabolic PDEs.
In the case of a general foliation, the topology of the leave through a point can change dramatically with the point; this gives many difficulties in studying the PRF. Therefore one may assume, at least at the first stage of study, the Riemannian manifold to be compact and to be fibred instead of being foliated.

Rovenski posed (in his project EU-FP7-P-2010-RG, No.~276919) the following conjecture:

\begin{conjecture} \label{RC}
Let a compact Riemannian manifold $(M^{n+p}, g)$ be the total
space of a smooth fiber bundle $\pi:M\rightarrow B$ with totally geodesic $p$-dimensional fibers.
Assume all mixed curvatures to be sufficiently close to a positive constant
(the degree of the closeness should be specified).
The PRF evolves the metric $g$, after the normalization, to a limit
metric whose mixed sectional curvature is a~function of a point, i.e., is independent of a plane.
\end{conjecture}

The conjecture seems to be an analogue of the following result by B\"{o}hm and Wilking.

\begin{theorem}[see Theorem 1.10 in \cite{ah}] \label{BWT}
On a compact manifold the Ricci flow evolves a Riemannian metric with 2-positive curvature operator
to a limit metric with constant sectional curvature.
\end{theorem}

Observe the following difference in statements of Conjecture \ref{RC} and Theorem \ref{BWT}: The sectional curvature of the limit metric is constant in Theorem \ref{BWT} while it can depend on a point in Conjecture \ref{RC}. The difference is caused by the absence of Schur's lemma in the case of fiber bundles.
Nevertheless, the statement of Conjecture \ref{RC} implies the Ferus inequality (\ref{FI}).

\subsection{One-dimensional foliations}
\label{sec:onedim}

In what follows, we will consider a one-dimensional foliation $\calf$ (i.e., $\calf$ is spanned by a unit vector field $N$) since this case is easier. Let $(M^{n+1}, g)$ be a~connected Riemannian manifold,
$N$ a~unit vector field on $TM$ and $\mathcal{D}$ its orthogonal codimension-one distribution.
In other words, $g(N,N)=1$ and $g(N,X)=0$ for $X\in\mathcal{D}$.
 In this case $p=1$, the~symmetric $(0,2)$-tensor in \eqref{E-rF} has the view
\begin{equation}\label{E-rN}
 r(X,Y) = R(X,N,Y,N)\qquad (X,Y\in TM).
\end{equation}
Its dual is the \textit{Jacobi operator} $R(N,\cdot)N$ for $N$.
 Now, $\mathcal{M}$ (called $N$-\textit{adapted metrics}) is the space of smooth Riemannian metrics
 on $M$ such that the vector field $N$ is unit.
The tensor \eqref{E-rN} provides an example of an $N$-\textit{truncated} symmetric $(0,2)$-tensor $S$, i.e.,
$S(N,X)=0$ for $X\in TM$.

For $p=1$, the PRF equation \eqref{E-GF-Rmix} says that for any vector $X$ orthogonal to $N$, the time derivative of $g(X,X)$ equals to minus twice the sectional curvature over the plane $X\wedge N$.

By Corollary~\ref{L-preserveN} (in Section~\ref{sec:tdep}), we have the following.

\begin{proposition}
The flow \eqref{E-GF-Rmix} preserves the length of $N$ and the orthogonality of vectors to $N$.
If~$N$ is a geodesic vector field for $t=0$ then it is geodesic for all $t$.
\end{proposition}

Notice that a circle bundle is a fiber bundle, where the fiber is a circle.

\begin{theorem}\label{T-main-loc}
Let $(M, g_0)$ be a closed Riemannian manifold with a unit vector field $N$.
Suppose that $N$-curves compose a circle bundle with  the total space $M$.
Then (\ref{E-GF-Rmix}) has a unique smooth solution $g_t$ defined on a positive time interval~$[0,\eps)$.
\end{theorem}

\section{Preliminaries}
\label{sec:prel}

In this section we survey the conullity and related tensors of extrinsic geometry of foliations,
describe their behavior under $N$-truncated variations of a metric,
and calculate the second derivative of the curvature tensor in the $N$-direction.

\subsection{Basic tensors of extrinsic geometry of a foliation}

The scalar \textit{second fundamental form} $h$ and the \textit{integrability tensor} $T$ of $\mathcal{D}$
are given~by
\begin{equation}\label{E-def-bT}
 h(X,Y)=\frac12\,g(\nabla_X Y+\nabla_Y X,\ N),\quad T(X,Y)=\frac12\,g([X,\,Y],\ N)\quad (X,Y\in\mathcal{D}).
\end{equation}
If $\mathcal{D}$ is integrable then $T=0$, and if $N$-curves compose a Riemannian foliation then $h=0$.
Notice that $g(\nabla_X Y,\ N)=h(X,Y)+T(X,Y)$.

The~(self-adjoint) \textit{Weingarten operator} $A: \mathcal{D}\to\mathcal{D}$ and the skew-symmetric operator $T^\sharp: \mathcal{D}\to\mathcal{D}$ are dual to $(0,2)$-tensors~$h$ and $T$, respectively:
\begin{equation}\label{E-shapeA}
 g(A(X),\,Y) = h(X,Y),\qquad g(T^\sharp(X),\,Y) = T(X,Y)\qquad (X,Y\in\mathcal{D}).
\end{equation}
The \textit{co-nullity tensor} $C: TM\to\mathcal{D}$ is defined~by
\begin{equation}\label{E-conulC}
  C(X)=-\nabla_{\!X} N,\qquad X\in TM.
\end{equation}
In particular, $C(N)=-\omega$, where $\omega=\nabla_{\!N} N$ is the \textit{curvature vector of $N$-curves}.

Let $^*$ be the conjugation of $(1,1)$-tensors on $\mathcal{D}$ with respect to $g$. Then
\begin{equation}\label{E-CC*}
 A=(C+C^{\,*})/2,\qquad T^\sharp =(C-C^{\,*})/2
\end{equation}
on $\mathcal{D}$.
Hence, $C=0$ when $\omega=0$
and $\mathcal{D}$ is integrable with totally geodesic integral manifolds
(i.e., $M$ is locally the direct product).

The self-adjoint
$(1,1)$-tensor $R_N$ given by $g(R_N(X),\,Y) = r(X,Y)$ is called the \textit{Jacobi operator} of $N$.
Its trace $\Ric_N =\tr\!R_N$ is the Ricci curvature in the $N$-direction.

\begin{lemma}\label{L-CC-riccati}
For a unit vector field $N$ $($i.e., $g(N,N)=1)$
we have, see for example $\cite[Lemma~2.25]{rov-m}$,
\begin{equation}\label{E-nablaT-1gen}
 \nabla_N\,C = C^{\,2} +R_N +g(\omega,\,\cdot)\,\omega - \nabla^\perp\omega.
\end{equation}
If $\omega=0$ then the symmetric and skew-symmetric parts of \eqref{E-nablaT-1gen} are, respectively,
 \begin{equation}\label{E-nablaT-1}
 \nabla_N\,A = A^2 +(T^\sharp)^2 +R_N,
 \qquad
 \nabla_N\,T^\sharp = A\,T^\sharp +T^\sharp A.
\end{equation}
\end{lemma}

\begin{remark}\label{R-contract}\rm
Tracing (\ref{E-nablaT-1gen})
and using $\Div\,N\!=\!-\tau_1$
yields the formula
\begin{equation}\label{E-IF1-RicN}
 N(\tau_1)=\tr (A^2) +\tr((T^\sharp)^2) +\Ric_N -\Div \omega,
\end{equation}
where
$\tau_1=\tr A=\tr_g h$ is the \textit{mean curvature} of $\mathcal{D}$.
Note that $\|T^\sharp\|^2=-\tr((T^\sharp)^2)$.
\end{remark}

\subsection{Time-dependent adapted metric}
\label{sec:tdep}

Let $S(g)$ be an $N$-\textit{truncated} symmetric $(0,2)$-tensor on $M$,
and $S^\sharp:TM\to TM$ its dual $(1,1)$-tensor. Consider
a family $g_t\in\mathcal{M}$ (with $0\le t<\eps$) of adapted metrics satisfying PDE
\begin{equation}\label{E-Sgeneral}
 \dt g=S(g).
\end{equation}
Since the difference of two connections is always a tensor, $\Pi_t:=\dt\nabla^t$ is a~$(1,2)$-tensor on $(M, g_t)$ with the symmetry $\Pi_t(X, Y)=\Pi_t(Y, X)$. It~is known, see~\cite{ah},
\begin{equation}\label{eq2}
 2\,g_t(\Pi_t(X, Y), Z)=(\nabla^t_X S_t)(Y,Z)+(\nabla^t_Y S_t)(X,Z)-(\nabla^t_Z S_t)(X,Y)
\end{equation}
for all $X,Y,Z\in\Gamma(TM)$.
If the vector fields $X=X(t),\,Y=Y(t)$ are $t$-dependent, then
\begin{equation}\label{eq2time}
 \dt\nabla^t_X Y = \Pi_t(X, Y) + \nabla_X(\dt Y)+ \nabla_{\dt X} Y.
\end{equation}

\begin{lemma}[see \cite{rw-m, rw4}]\label{L1-btAt2}
The tensors $A$, $T^\sharp$ and $C=A+T^\sharp$,
the mean curvature function $\tau_1$ of $\mathcal{D}$ and the curvature vector $\omega$ of $N$-curves evolve by (\ref{E-Sgeneral})~as
\begin{eqnarray}
\label{E1-S-A}
 \dt A\eq -\frac12\,\nabla_{N} S^\sharp +\frac12\,[A-T^\sharp,\ S^\sharp],\qquad
 \dt T^\sharp = -S^\sharp\,T^\sharp,\\
\label{E1-S-C}
 \dt C\eq -\frac12\,\nabla_N\,S^\sharp +\frac12\,[C,\ S^\sharp]-T^\sharp\,S^\sharp,\\
\label{E1-S-H}
 \dt\tau_1\eq -(1/2)\,N(\tr S^\sharp),\\
\label{E1-S-om}
 \dt\omega \eq -S^\sharp(\omega).
\end{eqnarray}
\end{lemma}

\noindent\textbf{Proof}.
For all $X,Y\in\mathcal{D}$, using (\ref{eq2}), we have
\begin{eqnarray*}
 2\,g\big(\dt(\nabla^t_X Y),\,{N}\big)
 \eq (\nabla^t_X S)(Y,{N})+(\nabla^t_Y S)(X,{N}) -(\nabla^t_{N} S)(X,Y)\\
 \eq -(\nabla^t_{N} S)(X,Y) -S(Y,\nabla^t_X {N}) -S(X,\nabla^t_Y{N}).
\end{eqnarray*}
From this and (\ref{E-def-bT}), we have the formula
\begin{equation}\label{E1-S-b}
 2\,\dt h(X,Y) = -(\nabla_N\,S)(X,Y)+S(C(X), Y)+S(X, C(Y))\qquad (X,Y\in\mathcal{D}).
\end{equation}
Note that $\dt T(X,Y)=0$.
Using identities (\ref{E-shapeA}), we find
\begin{equation*}
 g(\dt A(X),Y) = \dt h(X,Y) -(\dt g)(A(X),Y),\qquad
 g(\dt T^\sharp(X),Y) = -(\dt g)(T^\sharp(X),Y).
\end{equation*}
From the above, \eqref{E-Sgeneral} and (\ref{E1-S-b}) we deduce (\ref{E1-S-A}).
Using $\dt C=\dt A +\dt T^\sharp$ and (\ref{E1-S-A})$_1$ we obtain
\[
  \dt C= -\frac12\,\nabla_N\,S^\sharp +\frac12\,[A-T^\sharp,\ S^\sharp]-S^\sharp\,T^\sharp,
\]
and hence (\ref{E1-S-C}).
 Next, we deduce~(\ref{E1-S-H}):
 $\dt\tau_1=\tr(\dt A) = -\frac12\tr(\nabla_N\,S^\sharp)=-\frac12\,N(\tr S^\sharp)$.

From \eqref{eq2} we then have
 $2\,g_t(\dt\omega, X) = -2\,g_t(S^\sharp(\omega), X)$ for $X\in TM$
that is \eqref{E1-S-om}.
\qed

\begin{corollary}\label{L-preserveN}
The flow of metrics \eqref{E-Sgeneral} preserves the codimension-one distribution $\mathcal{D}$ orthogonal to $N$,
and $N$ is unit for all $g_t$. If $N$ is geodesic for $t=0$ then it is geodesic for all $t$.
\end{corollary}

\noindent\textbf{Proof}.
Since $S(N,\cdot)=0$, the flow (\ref{E-Sgeneral}) preserves the distribution $\mathcal{D}$
orthogonal to $N$, and $N$ is unit for all~$g_t$.
From \eqref{E1-S-om} and the theory of linear ODEs the last claim follows.
\qed

\begin{corollary}\label{C-ATC*}
For (\ref{E-Sgeneral}), the symmetries of $\dt A$, $\dt T^\sharp$ and $\dt R_N$ may be lost:
\begin{eqnarray}\label{E-dtAT*}
 && (\dt A)^* -\dt A = [S^\sharp,\ A],\qquad
 (\dt T^\sharp)^* +\dt T^\sharp=[T^\sharp, S^\sharp],\\
\label{E-dtRN*}
 && (\dt R_N)^* -\dt R_N = -[R_N, S^\sharp].
\end{eqnarray}
\end{corollary}

\noindent\textbf{Proof}. From (\ref{E1-S-A}) formulae (\ref{E-dtAT*})$_{1,2}$ follow.
Notice that
\begin{equation}\label{E-dtRdtr}
 \dt r(X,Y)= \dt(g(R_N(X),Y))
 =S(R_N(X),Y)+g((\dt R_N)(X),Y)\qquad (X,Y\in\mathcal{D}).
\end{equation}
From this and symmetry of $\dt\,r$ the equality \eqref{E-dtRN*} follows.
\qed


\vskip2mm
Denote $E$ the pull-back of $TM$ under the projection $M\times(0,\eps)\to M,\ (q,t)\to q$.
The fiber of $E$ over a point $(q,t)$ is given by $E_{(q,t)}=T_q M$.
A connection $\nabla$ on a vector bundle $E$ over $M$ is a~map
$\nabla:\mathcal{X}(M)\times\Gamma(E)\to\Gamma(E)$,
written as $(X, \sigma)\to\nabla_X \sigma$, that satisfies the following properties~\cite{ah}:

1. $\nabla$ is $C^\infty(M)$-linear in $X$: $\nabla_{f_1X_1+f_2X_2}\,\sigma=f_1\nabla_{X_1}\,\sigma+f_2\nabla_{X_2}\,\sigma$,

2. $\nabla$ is $\RR$-linear in $\sigma$:
$\nabla_X (\lambda_1\sigma_1+\lambda_2\sigma_2)=\lambda_1\nabla_{X}\,\sigma_1+\lambda_2\nabla_{X}\,\sigma_2$,
and

3. $\nabla$ satisfies the product rule:
$\nabla_X (f\sigma)=X(f)\,\sigma+f\nabla_{X}\,\sigma$.
\newline
A connection$\,\nabla$ on a vector bundle $E$ is said to be compatible
with a metric $g$ on $E$ if for any $\xi,\eta\in\Gamma(E)$ and $X\in\mathcal{X}(M)$,
we have $X(g(\xi,\,\eta))=g(\nabla_X\xi,\,\eta)+g(\xi,\,\nabla_X\eta)$.
Compatibility by itself is not enough to determine a unique connection.
 There is a natural connection $\widetilde\nabla$ on $E$, which extends the Levi-Civita connection on $TM$.
We need to specify only the covariant time derivative $\widetilde\nabla_{\dt}$.
Given any section $X$ of the vector bundle $E$, we define
\begin{equation}\label{E-tildenabla}
 \widetilde\nabla_{\dt} X=\dt X +\frac12\,S^\sharp(X)\ \ \mbox{for}\ \ X\in\mathcal{D},\qquad
 \widetilde\nabla_{\dt} N=0.
\end{equation}

\begin{lemma}
 The connection on $E$ is compatible with the natural bundle metric:
\begin{equation}\label{E-nablatg}
 \widetilde\nabla_{\dt}\, g =0.
\end{equation}
\end{lemma}

\noindent\textbf{Proof}.
 One may assume that $X,Y\in\mathcal{D}$ are constant in time. In this case, we have
 $\widetilde\nabla_{\dt} X=\frac12\,S^\sharp(X)$ and $\widetilde\nabla_{\dt} Y=\frac12\,S^\sharp(Y)$.
 Since $\dt g=S$, this and \eqref{E-Sgeneral} imply (\ref{E-nablatg}):
 \[
 (\widetilde\nabla_{\dt}\, g)(X,Y) = \dt g(X,Y) -g(\widetilde\nabla_{\dt} X, Y)-g(X, \widetilde\nabla_{\dt} Y)
 =(\dt g)(X,Y) -S(X,Y)=0.\qquad\qed
 \]

This connection is not symmetric: in general $\widetilde\nabla_{\dt} X\ne0$, while $\widetilde\nabla_{X}\dt =0$
always for $X\in\mathcal{D}$.
Clearly, the \textit{torsion tensor}
 ${\rm Tor}(X,Y):=\widetilde\nabla_XY-\widetilde\nabla_YX-[X,Y]$
vanishes if both arguments are spatial, so the only nonzero components are
\[
 {\rm Tor}(\dt, X) = \widetilde\nabla_{\dt} X -\widetilde\nabla_{X}\dt = \frac12\,S^\sharp(X)
 \qquad (X\in\mathcal{D}).
\]
However, each submanifold $M\times\{t\}$ is totally geodesic, so computing derivatives of spatial tangent
vector fields gives the same result as computing for sections of $T(M\times[0,\eps))$.
In particular, the corresponding Weingarten operators satisfy $\tilde A=A$.

\begin{remark}\rm
Using the connection (\ref{E-tildenabla}), we also have
\begin{eqnarray*}
 (\widetilde\nabla_{\dt} h)(X,Y) \eq \dt h(X,Y)-h(\widetilde\nabla_{\dt} X, Y)-h(X, \widetilde\nabla_{\dt} Y) = -\frac12(\nabla^t_{N} S)(X,Y),\\
 (\widetilde\nabla_{\dt} A)(X) \eq(\dt A)(X)-A(\widetilde\nabla_{\dt} X)
  = -\frac12(\nabla^t_{N} S^\sharp)(X)-\frac12\,[T^\sharp,\ S^\sharp].
\end{eqnarray*}
If $\mathcal{D}$ is integrable then, see (\ref{E1-S-A}),
\begin{equation*}
 \dt A = -\frac12\,\nabla_{N} S^\sharp +\frac12\,[A,\ S^\sharp],\qquad
 \widetilde\nabla_{\partial_t} A =-\frac12\,\widetilde\nabla_{N} S^\sharp.
\end{equation*}
\end{remark}

\subsection{The second $N$-derivative of the curvature tensor}

By Rectification Theorem, in a neighborhood of a point $q\in M$ there exist adapted coordinates
 $(x_0, x_1,\ldots,x_n)$ with coordinate fields $\partial_{x_i}=\partial_{i}$ for which $N$ is a coordinate field:
 $\partial_{x_0}=\partial_{0}=N$ and $g_{00}=1$, see \cite{cc1}.
The~Riemannian curvature tensor in components has the view $R_{ijkl}=g_{lp}R_{ijk}^p$,
where $R_{ijk}^p$ are components of the~curvature operator.

\vskip2mm
In analogy with \cite[Section~4.2.1]{ah}, define the quadratic in the curvature tensor $B\in \Lambda_0^4(M)$~as
\[
 B(X,Y,W,Z)=\<R(X,\ \cdot\ , Y,\,N),\   R(W,\ \cdot\ , Z,\,N)\>
 \quad\mbox{\rm for all}\quad X,Y,W,Z\in TM.
\]
In adapted coordinates (with respect to a unit vector field $N=\partial_0$), this becomes
 \[
  B_{ij kl} = g^{pq} R_{pi j0} R_{qk l0}.
 \]
Although generally $B_{ij kl}\neq B_{ji lk}$,
the tensor $B$ has some symmetries of the curvature tensor, as
\begin{equation}\label{E-symB}
 B_{ij kl} =B_{kl ij}.
\end{equation}

\begin{proposition}\label{P-SR-1}
 On a Riemannian manifold $(M,g)$ endowed with a unit vector field $N$,
 the second $N$-derivative of the curvature tensor satisfies
\begin{eqnarray}\label{E-rR0}
 \nabla^2_{N,N} R_{ijkl} \eq
 \nabla^2_{i,k}\,R_{j0l0}-\nabla^2_{j,k}\,R_{i0l0} +\nabla^2_{j,l}\,R_{i0k0} -\nabla^2_{i,l}\,R_{j0k0}\\
\nonumber
 \minus (B_{ij kl}-B_{ij lk}-B_{ji kl}+B_{ji lk}-2\,B_{kj li}+2\,B_{ki lj})
 +g^{pq}(R_{pj kl}\,r_{iq}-R_{pi kl}\,r_{jq})\,.
\end{eqnarray}
\end{proposition}

\noindent\textbf{Proof}.
Using the second Bianchi identity
$\nabla_{0} R_{ij\,kl} +\nabla_{i} R_{j0\,kl} +\nabla_{j} R_{0i\,kl}=0$ --
together with the linearity of $\nabla$ over the space of tensor fields -- we find that
\begin{equation}\label{E-SR-01}
 \nabla^2_{0,0} R_{ijkl}= \nabla^2_{0,i} R_{0jkl} -\nabla^2_{0,j} R_{0ikl}.
\end{equation}
It suffices to express first two terms in rhs of (\ref{E-SR-01}) using lower order terms.
To compute the first term in rhs of (\ref{E-SR-01}), we transpose $\nabla_0$ and $\nabla_i$,
\begin{equation}\label{E-SR-02}
 \nabla^2_{0,i} R_{0j kl}
 =\nabla^2_{i,0} R_{0j kl} +(R(\partial_i,\partial_0)R)(\partial_0,\partial_j,\partial_k,\partial_l).
\end{equation}
 We transform the first term in rhs of (\ref{E-SR-02}), using the second Bianchi identity
$\nabla_{0} R_{\,kl\,0j}+\nabla_{k} R_{\,l0\,0j}+\nabla_{l} R_{\,0k\,0j}=0$,
\begin{equation}\label{E-SR-03}
  \nabla^2_{i,0} R_{0j kl}= \nabla^2_{i,k} R_{j0 l0} -\nabla^2_{i,l} R_{0j 0k}.
\end{equation}
 Next, we transform the second term in the rhs of (\ref{E-SR-02}), using the identity
 $(R(X,Y)R)(Z,U,V,W)=-R(R(X,Y)Z,U,V,W)-R(Z,R(X,Y)U,V,W)-R(Z,U,R(X,Y)V,W)-R(Z,U,V, R(X,Y)W)$
 and noting that $R(X,Y)f=0$ where $f=R(Z,U,V,W)$,
 (in our case, $R(\partial_i,\partial_0)(R_{0j kl})=0$)
\begin{eqnarray}\label{E-SR-04}
\nonumber
 (R(\partial_i,\partial_0)R)(\partial_0,\partial_j,\partial_k,\partial_l) \eq -R(R(\partial_i,\partial_0)\partial_0,\partial_j,\partial_k,\partial_l)
 -\ldots -R(\partial_0,\partial_j,\partial_k, R(\partial_i,\partial_0)\partial_l)\\
\nonumber
 \eq -R_{i00}^{\,... q} R_{qj kl} -R_{i0j}^{\,... q} R_{0q kl} -R_{i0k}^{\,... q} R_{0j ql}
 -R_{i0l}^{\,... q} R_{0j kq}\\
 \eq g^{pq}(R_{i0 p0} R_{q j kl}+R_{0i jp} R_{0q kl} +R_{0i kp} R_{0j ql} +R_{0i lp} R_{0j kq}).
\end{eqnarray}
 The first term in the rhs of (\ref{E-SR-04}) is
\[
  g^{pq} R_{i0 p0} R_{qj kl}=-g^{pq} R_{jp kl}\,r_{iq}.  
\]
We transform the second term in the rhs of (\ref{E-SR-04}), using the first Bianchi identity,
\begin{eqnarray*}
  g^{pq} R_{0i jp} R_{0q kl} \eq - g^{pq} R_{pj i0} R_{kl q0}
  =g^{pq} R_{pj i0} R_{qk l0} +g^{pq} R_{pj i0} R_{lq k0} = B_{ji kl}-B_{ji lk}.
\end{eqnarray*}
The third and the fourth terms in the rhs of (\ref{E-SR-04}) are  transformed as
\begin{eqnarray*}
  g^{pq} R_{0i kp} R_{0j ql} +g^{pq} R_{0i lp} R_{0j kq}\eq
  -g^{pq} R_{pk i0} R_{ql j0} +g^{pq} R_{pl i0} R_{qk j0}
  =-B_{ki lj} +B_{li kj}.
\end{eqnarray*}
Hence (\ref{E-SR-04}) takes the following form:
 \begin{equation}\label{E-SR-05}
  (R(\partial_i,\partial_0)R)(\partial_0,\partial_j,\partial_k,\partial_l)
 = B_{ji kl}-B_{ji lk}-B_{ki lj}+B_{li kj} -g^{pq}R_{jp kl}\,r_{iq}.
\end{equation}
Substituting expressions of (\ref{E-SR-03}) and (\ref{E-SR-05}) into (\ref{E-SR-02}), we have
\begin{equation*}
 \nabla^2_{0,i} R_{0j kl} =\nabla^2_{i,k}\,R_{j0l0}-\nabla^2_{i,l}\,R_{j0k0}
 -(B_{ji lk} -B_{ji kl} +B_{ki lj}-B_{li kj}) -g^{pq}R_{jp kl}\,r_{iq}.
\end{equation*}
Using symmetry $i\leftrightarrow j$, we also have
\begin{eqnarray*}
 &&\nabla^2_{0,j} R_{0i kl} =\nabla^2_{j,k}\,R_{i0l0}-\nabla^2_{j,l}\,R_{i0k0}
 -(B_{ij lk} -B_{ij kl} +B_{kj li}-B_{lj ki}) -g^{pq}R_{ip kl}\,r_{jq}.
\end{eqnarray*}
By the above, (\ref{E-SR-01}) reduces to
\begin{eqnarray*}
 && \nabla^2_{0,0} R_{ijkl}
  =\nabla^2_{i,k}\,R_{j0 l0}-\nabla^2_{i,l}\,R_{j0 k0} -\nabla^2_{j,k}\,R_{i0 l0}+\nabla^2_{j,l}\,R_{i0 k0}
 +g^{pq}(R_{pj kl}\,r_{iq} -R_{pi kl}\,r_{jq})\\
 && -(B_{ji lk} -B_{ji kl} +B_{ki lj}-B_{li kj} -B_{ij lk} +B_{ij kl} -B_{kj li}+B_{lj ki}).
\end{eqnarray*}
Using the symmetry (\ref{E-symB}) of $B$, from the above we obtain \eqref{E-rR0}.
\qed

\begin{remark}\rm
The formula \eqref{E-rR0} becomes trivial for $j=l=0$.
\end{remark}

\section{Main results}

In this section we prove local existence/uniqueness theorem,
deduce the system of government equations for the curvature and conullity tensors
(which are parabolic along the leaves).

\subsection{Short-time existence and uniqueness}

To linearize the differential operator $g\to-2\,r(g)$, see (\ref{E-GF-Rmix}),
on the space of adapted metrics, we need the following.

\begin{proposition}[see \cite{ah}]\label{L-nablaNNR}
Let $g_t$ be a family of metrics on a manifold $M$ such that $\dt\,g =S$. Then
\begin{eqnarray}\label{E-dt RS}
\nonumber
 2\,\dt\,R(X,Y,Z,V) \eq \nabla^2_{X,V}S(Y,Z) +\nabla^2_{Y,Z}S(X,V)-\nabla^2_{X,Z}S(Y,V)-\nabla^2_{Y,V}S(X,Z)\\
\plus S(R(X,Y)Z,\ V)-S(R(X,Y)V,\ Z),
\end{eqnarray}
or in components,
\begin{eqnarray*}
 2\,\dt R_{ijkl}\eq \nabla^2_{i,l}\,S_{jk}+\nabla^2_{j,k}\,S_{il}-\nabla^2_{i,k}\,S_{jl}-\nabla^2_{j,l}\,S_{ik}
 +g^{pq}(R_{ijkp}S_{ql} -R_{ijlp}S_{qk}).
\end{eqnarray*}
\end{proposition}

Note that the first and second derivatives of a $(0,2)$-tensor $S$  are expressed~as
\begin{eqnarray}\label{E-12S}
\nonumber
 \nabla_{k}\,S_{jl}\eq \nabla_{k} S(\partial_j,\partial_l)
 =\partial_{k}(S(\partial_j,\partial_l)) -S(\nabla_k\partial_j,\partial_l)-S(\partial_j,\nabla_k\partial_l),\\
\nonumber
 \nabla^2_{i,k}\,S_{jl}\eq \nabla_{i}(\nabla_k\,S)_{jl} -\nabla_{\nabla_{i}\partial_k}\,S_{jl}\\
 \eq \nabla_{i}(\nabla_k\,S_{jl})
 -\nabla_k\,S(\nabla_{i}\partial_j, \partial_l)-\nabla_k\,S(\partial_j, \nabla_{i}\partial_l)
  -\nabla_{\nabla_{i}\partial_k}\,S_{jl}.
\end{eqnarray}

\begin{lemma}\label{L-dtRStrunc}
Let be $(M,g)$ a Riemannian manifold with a unit vector field $N$.
Then the tensor $r$ evolves by (\ref{E-Sgeneral}) with $N$-\textit{truncated} symmetric $(0,2)$-tensor $S(g)$
according to
\begin{eqnarray}\label{E-dt rS2}
\nonumber
 2\,\dt\,r_{ik} \eq -\nabla^2_{N,N}\,S_{ik}
 +\nabla_{N}\,S(C(\partial_i),\partial_k) +\nabla_{N}\,S(C(\partial_k), \partial_i)
 +S(C^2(\partial_i),\partial_k) +S(C^2(\partial_k), \partial_i) \\
\nonumber
 &&
 -2\,S(C(\partial_i),C(\partial_k)) +S(R_N(\partial_k), \partial_i) +S(R_N(\partial_i), \partial_k) \\
 &&
 -\nabla_{k}\,S(\omega, \partial_i) -\nabla_{i}\,S(\omega,\partial_k)
 -S(\nabla_i\,\omega,\partial_k)
 -S(\nabla_k\,\omega,\partial_i)
 +g(\omega,\partial_k)\,S(\omega, \partial_i),
\end{eqnarray}
where $\omega=\nabla_NN$.
If, in addition, $N$ is a geodesic vector field $($i.e., $\omega=0)$ then
\begin{eqnarray}\label{E-dtr-Sgeo}
\nonumber
 2\,\dt\,r_{ik} \eq -\nabla^2_{N,N}S_{ik}
 {+}\nabla_{N}S(C(\partial_i),\partial_k) {+}\nabla_{N}\,S(C(\partial_k), \partial_i)
 {+}S(C^2(\partial_i),\partial_k) {+}S(C^2(\partial_k), \partial_i)\\
 &&
 -2\,S(C(\partial_i),C(\partial_k)) +S(R_N(\partial_k), \partial_i) +S(R_N(\partial_i), \partial_k),\\
 \label{E-dtR-Sgeo}
  2\,\dt\,R_N \eq -\nabla^2_{N,N}\,S^\sharp
  {+}(\nabla_N S^\sharp) C {+}C^*\nabla_N S^\sharp {+}S^\sharp C^2 {+}(C^*)^2 S^\sharp
  {-}2\,C^* S^\sharp C {+}R_N S^\sharp {-}S^\sharp R_N.
\end{eqnarray}
\end{lemma}

\noindent\textbf{Proof}.
The symmetric $(0,2)$-tensor \eqref{E-rN} in adapted coordinates is $r(g)=(r_{ik}=R_{i 0 k 0})$.
We calculate the time derivative
 $\dt\,R_{i 0 k 0}=\dt\,r_{ik}$.
By Proposition~\ref{L-nablaNNR} with $j=l=0$, and using $R_{i00p}=-r_{ip},\ S_{q0}=0$, we then have
\begin{equation}\label{E-dt rS}
 2\,\dt\,r_{ik} = \nabla^2_{i,0}\,S_{0k}+\nabla^2_{0,k}\,S_{i0}-\nabla^2_{i,k}\,S_{00}-\nabla^2_{0,0}\,S_{ik}
 +g^{pq} r_{ip}\,S_{qk}.
\end{equation}
 By \eqref{E-conulC}, we have $\nabla_{0}\partial_i=\nabla_{i}\partial_0=-C(\partial_i)$.
By the above and \eqref{E-12S}, for a $N$-truncated symmetric $(0,2)$-tensor $S$, we have
 $\nabla_{k}\,S_{j0}= S(\partial_j,\, C(\partial_k))$, $\nabla_{k}\,S_{00}=0$
and
\begin{eqnarray*}
 \nabla^2_{i,0}\,S_{0k}\eq \nabla_{i}(\nabla_0\,S_{0k})
 -\nabla_0\,S(\nabla_{i}\partial_0, \partial_k)-\nabla_0\,S(\partial_0, \nabla_{i}\partial_k)
  -\nabla_{\nabla_{i}\partial_0}\,S_{0k}\\
 \eq\nabla_0\,S(C(\partial_i),\partial_k) +S(C^2(\partial_i),\partial_k)
 -\nabla_{i}\,S(\omega,\partial_k) -S(\nabla_i\,\omega,\partial_k),\\
 \nabla^2_{0,k}\,S_{i0}\eq \nabla_{0}(\nabla_k\,S_{i0})
 -\nabla_k\,S(\nabla_{0}\partial_i, \partial_0)-\nabla_k\,S(\partial_i, \omega)
 -\nabla_{\nabla_{0}\partial_k}\,S_{i0}\\
 \eq \nabla_0\,S(C(\partial_k), \partial_i) {-}\nabla_{k}\,S(\omega, \partial_i)
 {+}S(\nabla_0(C(\partial_k)), \partial_i)
 {+}S(C^2(\partial_k),\partial_i)
 \stackrel{\eqref{E-nablaT-1gen}}{=} \nabla_0\,S(C(\partial_k), \partial_i) \\
  \plus S(C^2(\partial_k), \partial_i) +S(R_N(\partial_k), \partial_i)
  +g(\omega,\partial_k)S(\omega, \partial_i) -\nabla_{k}\,S(\omega, \partial_i)
  -S(\nabla_k\,\omega,\partial_i),\\
 \nabla^2_{i,k}\,S_{00} \eq \nabla_{i}(\nabla_k\,S_{00})
 -2\,\nabla_k\,S(\nabla_{i}\partial_0, \partial_0) -\nabla_{\nabla_{i}\partial_k}\,S_{00}
 = 2\,S(C(\partial_i),C(\partial_k)).
\end{eqnarray*}
Hence \eqref{E-dt rS}, considered for a $N$-truncated tensor $S$, reduces to \eqref{E-dt rS2}.\qed

\vskip2mm\noindent\textbf{Proof of Theorem~\ref{T-main-loc}}.
We will use variations of the form $g(t)=g_0+t\,S$ with a $N$-truncated symmetric $(0,2)$-tensor $S$,
i.e., $S(N,\,\cdot)=0$.
We will show that $\nabla^2_{0,0}\,S_{ik}$ yields the principal symbol of order two,
and other terms are of order less then two.
By Lemma~\ref{L-dtRStrunc}, the linearization of $-2\,r$ is the second order differential operator
(elliptic along $N$-curves)
\[
 D(-2\,r)_{ik}= \nabla^2_{N,N}\,S_{ik} +\tilde Q_{ik},
\]
where $\tilde Q_{ik}$ consists of the first and zero order terms.
The~result then follows from the theory of parabolic PDEs on vector bundles,
see \cite[Section~5.1]{ah}, and the ``circle bundle" assumption.
\qed

\subsection{Evolution of the curvature tensor}

 If $M$ does not split along $N$,
the derivatives $\nabla_i\partial_j$ in adapted coordinates do not vanish simultaneously just at one point.
Indeed, $\nabla_0\partial_0=\omega=-C(\partial_0)$.

\begin{lemma}\label{L-NNS}
The difference of second derivatives,
$Q_{ik;jl}:=\nabla^2_{i,k}\,R_{j0l0}-\nabla^2_{i,k}\,r_{jl}$, has the view
\begin{eqnarray*}
 &&\hskip-20pt
 Q_{ik;jl}
 \!=\!\nabla_{i}R(\partial_j, C(\partial_k), \partial_l, \partial_0)
  {+}\nabla_{i}R(\partial_j, \partial_0, \partial_l, C(\partial_k))
  {+}\nabla_{k}R(\partial_j, C(\partial_i), \partial_l, \partial_0)
  {+}\nabla_{k}R(\partial_j, \partial_0, \partial_l, C(\partial_i))\\
 &&\hskip2pt
 {+}R(\partial_j, \!\nabla_{i} C(\partial_k), \partial_l, \partial_0)
 {+}R(\partial_j, \partial_0, \partial_l, \!\nabla_{i} C(\partial_k))
 {-}R(\partial_j, C(\partial_k), \partial_l, C(\partial_i))
 {-}R(\partial_j, C(\partial_i), \partial_l, C(\partial_k)).
\end{eqnarray*}
\end{lemma}

\noindent\textbf{Proof}. We use
\begin{eqnarray*}
 \nabla_{k}\,R_{j0l0} \eq\nabla_{k}\,r_{jl} -R(\partial_j, \nabla_{k}\partial_0, \partial_l, \partial_0)
  -R(\partial_j, \partial_0, \partial_l, \nabla_{k}\partial_0),\\
 \nabla^2_{i,k}\,r_{jl} \eq \nabla_{i}(\nabla_k\,r_{jl})
 -\nabla_k\,r(\partial_j, \nabla_{i}\partial_l)-\nabla_k\,r(\nabla_{i}\partial_j, \partial_l)
  -\nabla_{\nabla_{i}\partial_k}\,r_{jl}
\end{eqnarray*}
to calculate
\begin{eqnarray*}
 \nabla^2_{i,k}\,R_{j0l0} \eq \nabla_{i}(\nabla_{k}\,R)_{j0l0}-\nabla_{\nabla_{i}\partial_k}\,R_{j0l0}\\
 \eq\nabla_{i}\big(\nabla_{k}\,r_{jl} -R(\partial_j, \nabla_{k}\partial_0, \partial_l, \partial_0)
 -R(\partial_j, \partial_0, \partial_l, \nabla_{k}\partial_0)\big)
 -\nabla_{k}\,R(\nabla_{i}\partial_j, \partial_0, \partial_l, \partial_0) \\
 &&-\nabla_{k}\,R(\partial_j, \nabla_{i}\partial_0, \partial_l, \partial_0)
 -\nabla_{k}\,R(\partial_j, \partial_0, \nabla_{i}\partial_l, \partial_0)
 -\nabla_{k}\,R(\partial_j, \partial_0, \partial_l, \nabla_{i}\partial_0)\\
 &&-\nabla_{\nabla_{i}\partial_k}\,r_{j l}
 +R(\partial_j, \nabla_{\nabla_{i}\partial_k}\partial_0, \partial_l, \partial_0)
 +R(\partial_j, \partial_0, \partial_l, \nabla_{\nabla_{i}\partial_k}\partial_0)\\
 \eq \nabla_{i}\big(\nabla_{k}\,r_{jl}\big)
 -\nabla_{i}\,R(\partial_j, \nabla_{k}\partial_0, \partial_l, \partial_0)
 -\nabla_{i}\,R(\partial_j, \partial_0, \partial_l, \nabla_{k}\partial_0)
 -R(\nabla_{i}\partial_j, \nabla_{k}\partial_0, \partial_l, \partial_0)
 \\
 &&
 -R(\partial_j, \nabla_{i}\nabla_{k}\partial_0, \partial_l, \partial_0)
 -R(\partial_j, \nabla_{k}\partial_0, \nabla_{i}\partial_l, \partial_0)
 -R(\partial_j, \nabla_{k}\partial_0, \partial_l, \nabla_{i}\partial_0)\\
 &&-\nabla_{k}\,r(\nabla_{i}\partial_j, \partial_l)
 +R(\nabla_{i}\partial_j, \nabla_{k}\partial_0, \partial_l, \partial_0)
 +R(\nabla_{i}\partial_j, \partial_0, \partial_l, \nabla_{k}\partial_0) \\
 &&-\nabla_{k}\,r(\partial_j, \nabla_{i}\partial_l)
 +R(\partial_j, \nabla_{k}\partial_0, \nabla_{i}\partial_l, \partial_0)
 +R(\partial_j, \partial_0, \nabla_{i}\partial_l, \nabla_{k}\partial_0) \\
 &&-\nabla_{k}\,R(\partial_j, \nabla_{i}\partial_0, \partial_l, \partial_0)
  -\nabla_{k}\,R(\partial_j, \partial_0, \partial_l, \nabla_{i}\partial_0)
 -R(\nabla_{i}\partial_j, \partial_0, \partial_l, \nabla_{k}\partial_0)\\
 &&
 -R(\partial_j, \nabla_{i}\partial_0, \partial_l, \nabla_{k}\partial_0)
 -R(\partial_j, \partial_0, \nabla_{i}\partial_l, \nabla_{k}\partial_0)
 -R(\partial_j, \partial_0, \partial_l, \nabla_{i}\nabla_{k}\partial_0)\\
 &&-\nabla_{\nabla_{i}\partial_k}\,r_{j l}
 +R(\partial_j, \nabla_{\nabla_{i}\partial_k}\partial_0, \partial_l, \partial_0)
 +R(\partial_j, \partial_0, \partial_l, \nabla_{\nabla_{i}\partial_k}\partial_0)\\
 \eq \nabla^2_{i,k}\,r_{jl}
 +\nabla_{i}\,R(\partial_j, C(\partial_k), \partial_l, \partial_0)
 +\nabla_{i}\,R(\partial_j, \partial_0, \partial_l, C(\partial_k)) \\
 &&
  +\nabla_{k}\,R(\partial_j, C(\partial_i), \partial_l, \partial_0)
  +\nabla_{k}\,R(\partial_j, \partial_0, \partial_l, C(\partial_i))
   +R(\partial_j, \nabla_{i}\,C(\partial_k), \partial_l, \partial_0)\\
 &&
 -R(\partial_j, C(\partial_k), \partial_l, C(\partial_i))
 -R(\partial_j, C(\partial_i), \partial_l, C(\partial_k))
 +R(\partial_j, \partial_0, \partial_l, \nabla_{i}\,C(\partial_k)).
\end{eqnarray*}
The above yields the claim.
\qed

\begin{proposition}\label{T-RS-01} On a Riemannian manifold $(M,g)$ endowed with a unit vector field $N$,
the curvature tensor $R_{ijkl}\ (i,j,k,l\ge0)$ evolves by (\ref{E-GF-Rmix})
according to a heat type equation along $N$-curves
 \begin{eqnarray}\label{E-GFDDr-2}
  \dt\,R_{ijkl} \eq \nabla^2_{N,N}\,R_{ijkl}
  +B_{ij kl}-B_{ij lk}-B_{ji kl}+B_{ji lk}-2\,B_{kj li}+2\,B_{ki lj}\\
 \nonumber
  \minus g^{pq}(R_{qj kl}\,r_{ip}+R_{iq kl}\,r_{jp}+R_{ij ql}\,r_{kp}+R_{ij kq}\,r_{lp}) -\tilde Q\,,
\end{eqnarray}
where $\tilde Q=Q_{ik;jl}-Q_{il;jk}-Q_{jk;il}+Q_{jl;ik}$ $($with first order spatial derivatives, see Lemma~\ref{L-NNS}$)$.
 \end{proposition}

\noindent\textbf{Proof}. Applying Proposition~\ref{L-nablaNNR} with $S=-2\,r$, we have
\begin{equation}\label{E-compare1}
 \dt\,R_{ijkl} = \nabla^2_{i,k}\,r_{jl}-\nabla^2_{j,k}\,r_{il}+\nabla^2_{j,l}\,r_{ik}-\nabla^2_{i,l}\,r_{jk}
 -g^{pq}(R_{ij kp}\,r_{lq}+R_{ij pl}\,r_{kq}).
\end{equation}
Comparing (\ref{E-compare1}) with \eqref{E-rR0} completes the proof.
\qed

\begin{remark}\label{R-Q}
\rm By Lemma~\ref{L-NNS}, we find
\begin{eqnarray*}
 &&\hskip-6mm
 \tilde Q
  =\nabla_{C(\partial_k)}\,R_{ij0l}
  +\nabla_{C(\partial_j)}\,R_{i0 kl}
  +\nabla_{C(\partial_i)}\,R_{0j kl}
  +\nabla_{C(\partial_l)}\,R_{ij k0}\\
 &&\hskip-2mm
   + \nabla_{0}R(\partial_i, C(\partial_j), \partial_k, \partial_l)
  {+}\nabla_{0}R(C(\partial_i), \partial_j, \partial_k,\partial_l)
  {+}\nabla_{0}R(\partial_i, \partial_j, \partial_k, C(\partial_l))
  {+}\nabla_{0}R(\partial_i, \partial_j, C(\partial_k), \partial_l) \\
 &&\hskip-2mm
 {+}R(\partial_j, \!\nabla_{i} C(\partial_k), \partial_l, \partial_0)
 {+}R(\partial_j, \partial_0, \partial_l, \!\nabla_{i} C(\partial_k))
 {-}R(\partial_j, C(\partial_k), \partial_l, C(\partial_i))
 {-}R(\partial_j, C(\partial_i), \partial_l, C(\partial_k)) \\
 &&\hskip-2mm
 {-}R(\partial_j, \!\nabla_{i} C(\partial_l), \partial_k, \partial_0)
 {-}R(\partial_j, \partial_0, \partial_k, \!\nabla_{i} C(\partial_l))
 {+}R(\partial_j, C(\partial_l), \partial_k, C(\partial_i))
 {+}R(\partial_j, C(\partial_i), \partial_k, C(\partial_l)) \\
 &&\hskip-2mm
 {+}R(\partial_i, \!\nabla_{j} C(\partial_l), \partial_k, \partial_0)
 {+}R(\partial_i, \partial_0, \partial_k, \!\nabla_{j} C(\partial_l))
 {-}R(\partial_i, C(\partial_l), \partial_k, C(\partial_j))
 {-}R(\partial_i, C(\partial_j), \partial_k, C(\partial_l)) \\
 &&\hskip-2mm
 {-}R(\partial_i, \!\nabla_{j} C(\partial_k), \partial_l, \partial_0)
 {-}R(\partial_i, \partial_0, \partial_l, \!\nabla_{j} C(\partial_k))
 {+}R(\partial_i, C(\partial_k), \partial_l, C(\partial_j))
 {+}R(\partial_i, C(\partial_j), \partial_l, C(\partial_k)).
\end{eqnarray*}
\end{remark}

\begin{proposition} Let $(M,g)$ be a Riemannian manifold with a unit vector field $N$,
and $\omega=\nabla_N\,N$. Then the tensor $r$ evolves by (\ref{E-GF-Rmix}) according to
\begin{eqnarray}\label{E-dt-r_N}
\nonumber
 \dt\,r(X,Y) \eq \nabla^2_{N,N}\,r(X,\,Y) -\nabla_N\,r(C(X),Y) -\nabla_N\,r(X, C(Y)) \\
\nonumber
 \plus\nabla_X\,r(\omega, Y) +\nabla_Y\,r(\omega, X)
 +r(\nabla_X\,\omega, Y) +r(\nabla_Y\,\omega, X) -g(\omega, Y)\,r(X, \omega)\\
 \minus r(C^{\,2}(X), Y) - r(X, C^{\,2}(Y)) +2\,r(C(X),C(Y)) - 2\,r(X, R_N(Y)),
\end{eqnarray}
or in components,
 \begin{eqnarray}\label{E-GFDDric-2}
 \dt\,r_{ik} \eq\nabla^2_{N,N}\,r_{ik}
 -\nabla_{N}\,r(\partial_i, C(\partial_k)) -\nabla_{N}\,r(C(\partial_i), \partial_k) \\
\nonumber
 \minus r(\partial _i, C^2(\partial_k)) -r(C^2(\partial_i), \partial_k)
 +\,2\,r(C(\partial_i), C(\partial_k)) -2\,r(\partial_i, R_N(\partial_k))\\
\nonumber
 \plus\nabla_k\,r(\partial_i, \omega) +\nabla_i\,r(\omega, \partial_k) +r(\nabla_i\,\omega, \partial_k)
   -g(\omega, \partial_k)\,r(\omega, \partial_i).
\end{eqnarray}
If, in addition, $N$ is a geodesic vector field $($i.e., $\omega=0)$ then
\begin{eqnarray}\label{E-dt-r_N2}
 \dt\,r(X,Y)\eq\nabla^2_{N,N}\,r(X,Y) -\nabla_N\,r(C(X),Y) -\nabla_N\,r(X,C(Y))\\
\nonumber
 \minus r(C^{\,2}(X), Y) -r(X, C^{\,2}(Y)) +2\,r(C(X),C(Y)) -2\,r(X, R_N(Y)),
\end{eqnarray}
or in components,
\begin{eqnarray}\label{E-dtrRN-0}
\nonumber
 \dt\,r_{ik} \eq \nabla^2_{N,N}\,r_{ik}
 -\nabla_{N}\,r(C(\partial_i),\partial_k) -\nabla_{N}\,r(C(\partial_k), \partial_i)
 -r(C^2(\partial_i),\partial_k) -r(C^2(\partial_k), \partial_i) \\
 \plus 2\,r(C(\partial_i),C(\partial_k)) -2\,r(R_N(\partial_k), \partial_i),
\end{eqnarray}
furthermore, the Jacobi operator $R_N$ and its trace $\Ric_{N}$ evolve by (\ref{E-GF-Rmix}) according to
\begin{eqnarray}\label{E-dt-R_N2}
 \dt\,R_N \eq \nabla^2_{N,N}\,R_N-(\nabla_N R_N) C-C^{\,*}\nabla_N R_N -R_N C^{\,2}-(C^{\,*})^2 R_N+2\,C^{\,*}R_N C,\\
\label{E-dt-Ric_N2}
 \dt\Ric_N\eq N(N(\Ric_N)) -2\tr(A\,\nabla_N R_N) -4\tr((T^\sharp)^2\,R_N).
\end{eqnarray}
\end{proposition}

\noindent\textbf{Proof}.
We calculate the time derivative $\dt\,R_{i 0 k 0}=\dt\,r_{ik}$ and use
\eqref{E-GFDDr-2} with $j=l=0$,
\begin{eqnarray}\label{E-GFDDr-200a}
 \dt\,r_{i k}\eq\nabla^2_{N,N}\,r_{i k} +(B_{i0 k0}-B_{i0 0k}-B_{0i k0}+B_{0i 0k}-2\,B_{k0 0i}+2\,B_{ki 00})\\
\nonumber
 \minus g^{pq}(r_{kq}\,r_{ip}+R_{iq k0}\,r_{0p}+r_{i q}\,r_{kp}+R_{i0 kq}\,r_{0p}) +Q_{00;ik} -\tilde Q_{|\,j=l=0}\,.
\end{eqnarray}
Note that $  B_{i0 k0} = B_{i0 0k} = B_{0i k0} = B_{ki 00} = 0$ and $ B_{0i 0k} = g^{pq}\,r_{pi}\,r_{qk}$.
By Lemma~\ref{L-NNS}, we have
\begin{eqnarray*}
  Q_{ik,00} \eq -2\,r(C(\partial_i), C(\partial_k)),\\
  Q_{0k,i0} \eq -\nabla_0\,r(\partial_i, C(\partial_k)) -r(\partial_i, C^2(\partial_k))\\
  &&\hskip-1.5mm-\,r(\partial_i, R_N(\partial_k)) +\nabla_k\,r(\partial_i, \omega)
  -g(\omega, \partial_k)\,r(\omega, \partial_i) +r(\nabla_k\,\omega, \partial_i), \\
  Q_{i0,0k} \eq
  -\nabla_0\,r(C(\partial_i), \partial_k) -r(C^2(\partial_i), \partial_k)
  +\nabla_i\,r(\partial_k, \omega) +r(\nabla_i\,\omega, \partial_k).
\end{eqnarray*}
(Note that $Q_{00; ik}=0$ when $\omega=0$).
Hence
\begin{eqnarray*}
 &&\hskip-10mm Q_{00;ik}-\tilde Q_{|\,j=l=0} = Q_{i0;0k} + Q_{0k;i0} -Q_{ik,00} \\
 && = -\nabla_0\,r(\partial_i, C(\partial_k)) -\nabla_0\,r(C(\partial_i), \partial_k)
    -r(\partial_i, C^2(\partial_k)) -r(C^2(\partial_i), \partial_k)  -r(\partial_i, R_N(\partial_k)) \\
 && +2\,r(C(\partial_i), C(\partial_k))
 +\nabla_k\,r(\partial_i, \omega) +\nabla_i\,r(\omega, \partial_k) +r(\nabla_i\,\omega, \partial_k)
    +r(\nabla_k\,\omega, \partial_i) -g(\omega, \partial_k)\,r(\omega, \partial_i).
\end{eqnarray*}
Thus, \eqref{E-GFDDr-200a} reads as (\ref{E-GFDDric-2}).
The formula
\eqref{E-dt-r_N} with $\omega=0$ is equivalent to (\ref{E-dt-r_N2}).
By \eqref{E-dt-r_N2} and \eqref{E-dtRdtr} for $S=-2\,r$ we get \eqref{E-dt-R_N2}.
Tracing \eqref{E-dt-R_N2} and using $\dt(\tr R_N)=\tr(\dt R_N)$
and $\tr(T^\sharp\nabla_N\,R_N)=0$,~yields
\begin{equation*}
 \dt\Ric_N = N(N(\Ric_N)) -2\tr(C\,\nabla_N\,R_N) -4\tr((C\,T^\sharp\,R_N).
\end{equation*}
From this and the property $\tr(B_1B_2)=\tr(B_2B_1)$ we have \eqref{E-dt-Ric_N2}.
\qed

\begin{remark}\rm
Alternatively, one may deduce \eqref{E-GFDDric-2} substituting  $S=-2\,r$ into \eqref{E-dt rS2} of Lemma~\ref{L-dtRStrunc}.
\end{remark}

\noindent
We apply Uhlenbeck's trick (see \cite{ah}) to remove a group of terms in \eqref{E-GFDDr-2} with a 'change of variables'.

\begin{corollary}
If the metric $g_{ij}$ evolves by \eqref{E-GF-Rmix} then the curvature $R_{ijkl}$ evolves according to
 \begin{equation}\label{E-nablatRR}
  \widetilde\nabla_{t}\,R_{ijkl} = \nabla^2_{N,N}\,R_{ijkl}
  +B_{ij kl}-B_{ij lk}-B_{ji kl}+B_{ji lk}-2\,B_{kj li}+2\,B_{ki lj} -\tilde Q
\end{equation}
(where $\tilde Q$ is given in Proposition~\ref{T-RS-01}),
and the tensor $r_{ik}$ evolves according to
\begin{eqnarray}\label{E-nablatrr}
\nonumber
 \widetilde\nabla_{t}\,r_{ik} \eq \nabla^2_{N,N}\,r_{ik}
 -\nabla_{N}\,r(C(\partial_i),\partial_k) -\nabla_{N}\,r(C(\partial_k), \partial_i)\\
 \minus r(C^2(\partial_i),\partial_k) -r(C^2(\partial_k), \partial_i) +2\,r(C(\partial_i),C(\partial_k))\\
\nonumber
 \plus\nabla_k\,r(\partial_i, \omega) +\nabla_i\,r(\omega, \partial_k) +r(\nabla_i\,\omega, \partial_k)
   +r(\nabla_k\,\omega, \partial_i) -g(\omega, \partial_k)\,r(\omega, \partial_i).
\end{eqnarray}
\end{corollary}

\noindent\textbf{Proof}.
Using definition
$\widetilde\nabla_{\partial_t}\partial_i= \dt(\partial_i) -R_N(\partial_i)=-R_N(\partial_i)$,
we obtain
\begin{eqnarray*}
 \widetilde\nabla_{t}\,R_{ijkl}
 \eq\partial_t\,R_{ijkl}
 -R(-\!\widetilde\nabla_{\partial_t}\partial_i, \partial_j, \partial_k, \partial_l)
 -\ldots
 -R(\partial_i, \partial_j, \partial_k, -\!\widetilde\nabla_{\partial_t}\partial_l)\\
 \eq
 \partial_t\,R_{ijkl} +g^{pq}(r_{iq}R_{pjkl} +r_{jq}R_{ipkl}+r_{kq}R_{ijpl}+r_{lq}R_{ijkp}).
\end{eqnarray*}
From this and \eqref{E-GFDDr-2} the equation \eqref{E-nablatRR} follows.
Similarly, \eqref{E-nablatrr} follows from \eqref{E-GFDDric-2}.
\qed

\subsection{Evolution of the co-nullity tensor}

For (\ref{E-GF-Rmix}) with a unit vector field $N$,
using Corollary~\ref{C-ATC*} with $S=-2\,r$, we obtain
\begin{equation}\label{EGF-dtAT*}
 (\dt A)^* -\dt A = 2\,[A, \,R_N],\qquad
 (\dt T^\sharp)^* +\dt T^\sharp=-2\,[T^\sharp, \,R_N].
\end{equation}
We also have $(\dt R_N)^*=\dt R_N$.

\begin{proposition}\label{L1-btAt3}
For a unit geodesic vector field $N$, the tensors $A$, $T^\sharp$ and $C$
and the mean curvature function $\tau_1$ of $\mathcal{D}$
evolve by (\ref{E-GF-Rmix}) according to
\begin{eqnarray}
\label{E1-b-A}
 \dt A\eq\nabla_{N}(\nabla_{N}A) -2\,A\,\nabla_{N} A +[A^2,\ T^\sharp] -2(T^\sharp)^2 A -2\,T^\sharp A T^\sharp
 ,\\
\label{E1-b-T}
 \dt T^\sharp \eq 2\,(\nabla_N\,A) T^\sharp -2\,A^2 T^\sharp -2\,(T^\sharp)^3,\\
\label{E1-b-C}
\dt C\eq \nabla_{N}(\nabla_{N} C) -(C+C^{\,*})\,\nabla_{N} C -(C-C^{\,*}) C^{\,2}
 ,\\
\label{E1-b-H}
 \dt \tau_1 \eq N(N(\tau_1))-N(\tr (A^2))-4\tr(A\,T^\sharp)^2.
\end{eqnarray}
\end{proposition}

\noindent\textbf{Proof}. From \eqref{E1-S-A}--\eqref{E1-S-C}
 (of Lemma~\ref{L1-btAt2} with $S=-2\,r$) we obtain
\begin{eqnarray}\label{E-rN-A}
 \dt A\eq \nabla_N\,R_N +[T^\sharp -A,\ R_N],\qquad
 \dt T^\sharp = 2\,R_N T^\sharp,\\
\label{E-rN-C}
 \dt C\eq \nabla_N\,R_N +[R_N, C] +2\,T^\sharp R_N.
\end{eqnarray}
Substituting $R_N$ from \eqref{E-nablaT-1gen} into (\ref{E-rN-A})--\eqref{E-rN-C}, we obtain (\ref{E1-b-A})--(\ref{E1-b-C}).
Then, tracing \eqref{E1-b-A} yields \eqref{E1-b-H}.
\qed

\begin{remark}\rm
From \eqref{E1-S-H} with $S=-2\,r$ we have
\begin{equation}\label{E-rN-tau1}
 \dt\tau_1 = N(\Ric_N).
\end{equation}
Substituting $\Ric_N$ from \eqref{E-IF1-RicN} with $\omega=0$ into \eqref{E-rN-tau1} yields
\begin{equation}\label{E1-b-H-one}
 \dt \tau_1 = N(N(\tau_1))-N(\tr (A^2))-N(\tr((T^\sharp)^2).
\end{equation}
One may use \eqref{E-nablaT-1}$_2$ to calculate
\[
 \nabla_N((T^\sharp)^2)=A(T^\sharp)^2+(T^\sharp)^2 A+2\,T^\sharp A T^\sharp
 \ \Longrightarrow\ N(\tr((T^\sharp)^2))=4\tr(A(T^\sharp)^2).
\]
The above again yields \eqref{E1-b-H}.
\end{remark}

\begin{corollary}\label{P-nablaA}
Let $N$ be a unit geodesic vector field with integrable orthogonal distribution.
Then the Weingarten operator $A$ evolves by \eqref{E-GF-Rmix} according to
\begin{equation}\label{Ef-AE}
 \widetilde\nabla_{\dt} A = \nabla_N(\nabla_N\,A) -\nabla_N (A^2).
\end{equation}
\end{corollary}

\noindent\textbf{Proof}.
Due to definition \eqref{E-tildenabla}, we obtain $\widetilde\nabla_{\partial_t}= \dt -R_N$.
If $\mathcal{D}$ is integrable then (\ref{E1-b-A}) reads as
\begin{equation}\label{Ef-AE1}
 \dt A =\nabla_{N}(\nabla_{N} A) -2\,A\,\nabla_{N} A.
\end{equation}
One may assume that $X\in\mathcal{D}$ is constant in time. In this case, replacing $R_N$ due to (\ref{E-nablaT-1})$_1$, we have
\begin{eqnarray*}
 (\widetilde\nabla_{\dt} A)(X) \eq \widetilde\nabla_{\dt}(A(X)) - A(\widetilde\nabla_{\dt} X)
 =\dt(A(X))-R_N A (X) - A(-R_N (X)) \\
 \eq \big( \dt A-(\nabla_N\,A) A + A\,\nabla_N\,A\big)(X)
\end{eqnarray*}
for any $X\in\mathcal{D}$.
Applying (\ref{Ef-AE1}), we obtain (\ref{Ef-AE}).
\qed

\section{Examples}


In this section we show that PRF preserve several classes of foliations
and prove convergence of solution metrics under certain conditions.
 Let $N$ be a geodesic unit vector field and let $\Phi$ be a~function on $M$ satisfying $N(\Phi)=0$.
For \textit{normalized PRF} \eqref{E-GF-Rmix-Phi}, the co-nullity and the integrability tensors
evolve, see \eqref{E1-b-T}, \eqref{E1-b-C}, according~to
\begin{eqnarray}\label{E1-b-TPhi}
 \dt T^\sharp \eq 2\,(\nabla_N\,A) T^\sharp -2\,A^2 T^\sharp -2\,(T^\sharp)^3 -2\,\Phi\,T^\sharp,\\
\label{E1-b-CPhi}
 \dt C \eq \nabla_{N}(\nabla_{N} C) -(C+C^{\,*})\,\nabla_{N} C -(C-C^{\,*}) C^{\,2}
 -2\,\Phi\,T^\sharp,
\end{eqnarray}
while the Weingarten operator $A$ evolves by \eqref{E1-b-A}.

\begin{lemma}
Let $g(t),\ t\in[0,T]$ be a solution to \eqref{E-GF-Rmix}. If $|R_N|\le C$ on $M\times[0,T]$ then
\[
 e^{-2C\,t}g(0)_{|\,\mathcal D} \le g(t)_{|\,\mathcal D} \le e^{2C\,t}g(0)_{|\,\mathcal D}\quad
 \mbox{\rm for all} \ t\in[0,T].
\]
\end{lemma}

\noindent\textbf{Proof}.
This is similar to one of \cite[Lemma~8.5]{ah}.
\qed

\subsection{Evolution of a Riemannian geodesic foliation}

By the existence/uniqueness Theorem~\ref{T-main-loc}, \eqref{E1-S-om} and \eqref{E1-b-A}, we have the following.

\begin{proposition}
If  $\omega=0$ and $A=0$ at $t=0$ then the flow \eqref{E-GF-Rmix-Phi} preserves these properties.
\end{proposition}

 Suppose that $N$-curves compose a geodesic Riemannian foliation and $T\ne0$.
(Examples of such foliations are Hopf fibrations of odd-dimensional spheres).
 By Lemma~\ref{L-CC-riccati}, we then have
\[
 \nabla_N\,T^\sharp=0,\qquad R_N=-(T^\sharp)^2\ge0.
\]
This yields $\nabla_NR_N=0$, $\nabla_N\,r=0$ and $N(\Ric_N)=0$, hence \eqref{E-GF-Rmix-Phi} reduces to ODE in the variable $t$.

\begin{theorem}
Let $(M, g_0)$ be a Riemannian manifold with a unit geodesic vector field $N$ such that $A=0$ and $T\ne0$.
 If $r\le\Phi\,\hat g$ and $r_{|\,\mathcal D}>0$ at $t=0$ then \eqref{E-GF-Rmix-Phi} admits a~unique global solution $g_t\ (t\in\RR)$ such that $\lim\limits_{t\to-\infty} R_N(t)=\Phi\,\hat\id$ and $\lim\limits_{t\to\infty}R_N(t)=0$.
\end{theorem}

\noindent\textbf{Proof}.
By \eqref{E-dtR-Sgeo} with $S^\sharp=-2\,R_N +2\,\Phi\,\hat\id$, $C=T^\sharp$ and $\omega=0$, we obtain
\begin{eqnarray*}
 \dt\,R_N \eq -R_N (T^\sharp)^{\,2} -(T^\sharp)^2 R_N -2\,T^\sharp R_N T^\sharp -4\,\Phi R_N
 =4\,R_N(R_N -\Phi\,\hat\id),\\
 \dt\Ric_N \eq -4\tr((T^\sharp)^2 R_N) -4\,\Phi\Ric_N
 =4\tr(R_N^2) -4\,\Phi\Ric_N \ge \frac4n\,(\Ric_N)^2 -4\,\Phi\Ric_N.
\end{eqnarray*}
One may show  that \eqref{E-GF-Rmix-Phi} preserves the positive $\Ric_N$.
By Proposition~\ref{L1-btAt3}, we also have
\begin{equation}\label{E1-b-T-Riem}
 \dt T^\sharp = -2\,T^\sharp\big((T^\sharp)^2 +\Phi\,\hat\id\big).
\end{equation}
In our case $r_{|\,\mathcal D}>0$, the dimension $n$ should be even
(indeed, if $n$ is odd then the skew symmetric operator $T^\sharp$ has zero eigenvalues,
hence $R_N$ also has zero eigenvalues).

Let $\mu_i(t)>0$ be the eigenvalue and $e_i(t)$ the eigenvector of $R_N(t)$ under the flow~\eqref{E-GF-Rmix-Phi}. Then
\[
 \dt e_i = (\mu_i - \Phi) e_i,\qquad
 \dt \mu_i = 4\mu_i(\mu_i - \Phi).
\]
Hence the PRF preserves the directions of eigenvectors of $R_N$.
Furthermore,
if $\Phi\ge\mu_i(0)>0$ then $\lim\limits_{t\to-\infty}\mu_i(t)=\Phi$ and $\lim\limits_{t\to\infty}\mu_i(t)=0$.
\qed

\subsection{Evolution of warped product metrics}
\label{sec:rot_sym}

Let us look at what happens for a general \textit{warped product metric}
$g={\rm dx}^2+\varphi^2(x)\bar g$ on $M=[0,l]\times\bar M$, where $(\bar M^n,\bar g)$ is a Riemannian manifold
and $l$ a positive real, see \cite{Pet}.
(The rotational symmetric metrics, i.e., $\bar M$ is a unit $n$-sphere, are the particular case;
such metrics appear on surfaces of revolution in space forms).
The submanifolds $\{x\}\times\bar M$ compose a totally umbilical foliation on $M$ with a unit normal $N=\partial_x$. We have
\begin{eqnarray*}
 \omega \eq 0,\qquad
 T=0,\qquad
 A=-(\varphi_{,x}/\varphi)\,\hat\id,\\
 r \eq -(\varphi_{,xx}/\varphi)\,\hat g,\quad
 R_N = -(\varphi_{,xx}/\varphi)\,\hat\id,\quad
 \Ric(N,N)=-n\,\varphi_{,xx}/\varphi\quad({\rm when}\ \ \varphi\ne0).
\end{eqnarray*}
Thus, $K(N,X)=-\varphi_{,xx}/\varphi$ for $X\perp N$.
We apply the existence/uniqueness Theorem~\ref{T-main-loc}
to conclude that the flow \eqref{E-GF-Rmix-Phi} preserves totally umbilical foliations (with $N$ the unit normal).

\begin{proposition}\label{L-C0}
The flow \eqref{E-GF-Rmix-Phi} preserves each of the properties $C=\lambda\,\hat\id$ and $C=0$.
\end{proposition}

Thus we have the following.

\begin{corollary}
The flow \eqref{E-GF-Rmix-Phi} preserves warped product metrics.
\end{corollary}

Now, let a family of warped product metrics $g_t={\rm dx}^2+\varphi^2(t,x)\bar g$
solves \eqref{E-GF-Rmix-Phi} on $M$.
This yields the boundary value problem for the warping function $\varphi$,
\begin{equation}\label{E-r}
 \dt\varphi=\varphi_{,xx}+\Phi\,\varphi,\quad \varphi(0,x)=\varphi_0(x),\quad
 \varphi(t,0)=\mu_0(t),\quad \varphi(t,l)=\mu_1(t),
\end{equation}
where $\mu_j(t)\ge 0\ (j=0,1)$. The Cauchy's problem \eqref{E-r} has a
unique classical solution $\varphi(t,\,\cdot)$ for all $t\ge0$.
We are interested in convergence of solutions of this problem to a stationary state.

 Assume that the functions $\mu_j(t)\;(j=0,1)$ are continuously differentiable on $[0,\infty)$,
 and there exist limits
$\lim_{\,t\rightarrow\infty}\mu_j(t)=\tilde\mu_j$ and $\lim_{\,t\rightarrow\infty}\mu_j^\prime(t)=0$.

\begin{lemma}\label{Ex-tilder}
Let $\tilde\varphi(x)$ be a solution of the stationary Cauchy's problem on $[0,\,l]$
\begin{equation}\label{statprob}
 \tilde\varphi_{\,,xx}+\Phi\,\tilde\varphi=0,\qquad
 \tilde\varphi(0)=\tilde\mu_0,\qquad
 \tilde\varphi(l)=\tilde\mu_1.
\end{equation}
Then
\[
 \tilde\varphi(x)=\left\{\begin{array}{cc}
  \frac{\tilde\mu_1\sin(\sqrt{\Phi}\,x) +\,\tilde\mu_0\sin(\sqrt{\Phi}\,(l-x))}{\sin(\sqrt{\Phi}\,l)} &
 {\rm if} \ \ 0<\Phi<(\pi/l)^2,\\
 \tilde\mu_0+(\tilde\mu_1-\tilde\mu_0)(x/l) & \quad {\rm if}\ \Phi=0,\\
  \frac{\tilde\mu_1\sinh(\sqrt{-\Phi}\,x) +\,\tilde\mu_0\sinh(\sqrt{-\Phi}\,(l-x))}{\sinh(\sqrt{-\Phi}\,l)} &
 \quad {\rm if}\ \Phi<0.\\
 \end{array}\right.
\]
For the resonance case of $\Phi=(\pi/l)^2>0$, the problem \eqref{statprob} is solvable if and only if $\tilde\mu_1=-\tilde\mu_0$, and in this case the solutions are
$\tilde\varphi(x)=C\sin(\pi x/l)+\tilde\mu_0\cos(\pi x/l)$ where $C>0$ is constant.
\end{lemma}

 Denote
\begin{equation*}
  \delta_j(t):=\mu_j(t)-\tilde\mu_j\;(j=0,1),\qquad
  U(t,x):=\delta_0(t)\frac{l-x}{l}+\delta_1(t)\frac{x}{l},
\end{equation*}
\begin{equation}\label{dfv0U}
 v_0(x):=\varphi_0(x)-\tilde\varphi(x)-U(0,x),
\end{equation}
\begin{equation}\label{dfftx}
 f(t,x):=\Phi U -\dt U
 = (\Phi\delta_0-\delta'_0) +\frac xl\big( \Phi(\delta_1-\delta_0) +\delta'_0-\delta'_1\big).
\end{equation}
Consider the Fourier series
\begin{equation*}
 f(t,x)=\sum\nolimits_{j=1}^\infty f_j(t)\sin(\pi j\,x/l),\quad
 v_0(x)=\sum\nolimits_{j=1}^\infty v_j^0\sin(\pi j\,x/l),
\end{equation*}
where $f_j(t)=\frac{2}{l}\int_0^lf(t,s)\sin(\pi j s/l)\ds$
and $v_j^0=\frac{2}{l}\int_0^lv_0(s)\sin(\pi j s/l)\ds$.

For the warped product initial metric, the global solution of \eqref{E-GF-Rmix-Phi} converges to one with constant mixed sectional curvature.

\begin{theorem}\label{T-stabnonstat}
Let the warped product metrics $g_t$ on $M=[0,l]\times\bar M$ solve \eqref{E-GF-Rmix-Phi}.
If~$\Phi>(\pi/l)^2$ then the metrics diverge as $t\to\infty$, otherwise they converge uniformly for $x\in[0,l]$ to the limit metric $g_\infty$ with $r(g_\infty)=\Phi\,\hat g_\infty$,
hence, the mixed sectional curvature of $g_\infty$ is constant.
Certainly,

$(i)$ if $\Phi<(\pi/l)^2$ then $g_\infty =
{\rm dx}^2+\tilde\varphi^2(x)\,\bar g$, see Lemma~\ref{Ex-tilder}.

$(ii)$ if $\Phi=(\pi/l)^2$ and additional assumptions hold
\begin{equation}\label{addcond}
 \tilde\mu_j=0,\qquad
 \int_0^\infty\big(|\delta_j(\tau)|
 +|\delta_j^{\,\prime}(\tau)|\big)\dtau<\infty\qquad (j=0,1),
\end{equation}

then $g_\infty ={\rm dx}^2+\varphi_\infty^2(x)\bar g$
and $\varphi_\infty=\big(v_1^0+\int_0^\infty f_1(\tau)\dtau\big)\sin(\pi x/l)$.
\end{theorem}

\noindent\textbf{Proof}.
Let $\varphi(t,x)$ be a global solution of (\ref{E-r}).
Note that $\nu(t)\to0$ as $t\to\infty$, where
\[
 \nu(t):=|\Phi|(|\delta_0(t)|+|\delta_1(t)|)+|\delta_0^{\,\prime}(t)|+|\delta_1^{\,\prime}(t)|.
\]

$(i)$ Denote
\begin{eqnarray*}
 M_0(t) \eq\Big(\sum\nolimits_{j=1}^\infty e^{\,2\,(\pi/l)^2(1-j^2)\,t}\Big)^{1/2},\quad
 M_1(t) = \frac{6}{\pi}\sum\nolimits_{j=1}^\infty e^{\,(1-\theta)\,(\pi/l)^2(1-j^2)\,t}\,, \\
 M_2 \eq\frac{6}{\pi}\sum\nolimits_{j=1}^\infty \frac1{j\,((\pi j/l)^2-\Phi)}\,.
\end{eqnarray*}
We shall prove the following estimate which implies the claim~$(i)$:
\begin{eqnarray}\label{estdifference}
 &&\hskip-12mm|\varphi(t,x)-\tilde\varphi(x)|\le \max\{|\delta_0(t)|,|\delta_1(t)|\}
 + M_0(t) e^{\,(\Phi-(\pi/l)^2)\,t}\,\Vert v_0\Vert_{L_2}\\
\nonumber
 && +\, M_1(t)\big((\pi/l)^2-\Phi\big)^{-1} e^{(1-\theta)\,(\Phi-(\pi/l)^2)\,t}\sup\nolimits_{\,\tau\in[0,\,\theta t]}\nu(\tau)
 + M_2\sup\nolimits_{\,\tau\in[\theta t,\,t]}\nu(\tau)
\end{eqnarray}
for any $\theta\in(0,1)$ and all $(t,x)\in[0,\infty)\times[0,l]$.
 It is easy to check that the function
\begin{equation}\label{dfvtx}
 v(t,x)=\varphi(t,x)-\tilde\varphi(x)-U(t,x)
\end{equation}
is the solution of the following Cauchy's problem:
\begin{equation}\label{boundprv}
\partial_tv=v_{\,,xx}+\Phi v+f(t,x),\quad v(0,x)=v_0(x),\quad v(t,0)=v(t,l)=0.
\end{equation}
We may write
 $v(t,x)=\sum\nolimits_{j=1}^\infty v_j(t)\sin(\pi j\,x/l)$.
Substitution the series into (\ref{boundprv})$_{1,2}$ and comparison of the
coefficients of the series yield the Cauchy's problem for $v_j(t)$:
\begin{equation}\label{Cauchvj}
 v_j^{\,\prime}=(\Phi-(\pi j/l)^2)v_j+f_j(t),\qquad v_j(0)=v^0_j.
\end{equation}
Using series
 $1=\sum\nolimits_{j\ge1}\frac{2}{\pi j}\,(1-(-1)^j)\sin(\frac{\pi j\,x}l)$ and
 $x=-\sum\nolimits_{j\ge1}\frac{2 l}{\pi j}\,(-1)^{j}\sin(\frac{\pi j\,x}l)$,
from \eqref{dfftx} we find
 $f_j(t) = \frac{2}{\pi j}\,\big((1-(-1)^j)(\Phi\delta_0-\delta'_0)
 -(-1)^{j}\big(\Phi(\delta_1-\delta_0)+\delta'_0-\delta'_1\big)\big)$.
Hence,
\begin{equation}\label{E-fseries}
 |f_j(t)|
 \le\frac{2}{\pi j}\,\big(2\,(|\Phi\delta_0|+|\delta'_0|)
 +|\Phi|(|\delta_1|+|\delta_0|)+|\delta'_0|+|\delta'_1|\big)
 \le\frac{2}{\pi j}\,\big(2\,\nu(t)+\nu(t)\big) \le \frac{6}{\pi j}\,\nu(t).
\end{equation}
By Lemma~\ref{L-stabODE} with $a=\Phi-(\pi j/l)^2$ and $\nu(t)=f_j(t)$,
for any $\theta\in(0,1)$ we get the estimate
\begin{equation}\label{estvj}
 |v_j(t)|\le|v_j^0| e^{\,(\Phi-(\pi j/l)^2)\,t} -(\Phi-(\pi j/l)^2)^{-1}\big(
 e^{(1-\theta)(\Phi-(\pi j/l)^2)\,t}\!\sup\limits_{\,\tau\in[0,\,\theta t]}|f_j(\tau)|
 +\!\sup\limits_{\,\tau\in[\theta t,\,t]}|f_j(\tau)|\big).
\end{equation}
The above and \eqref{E-fseries}\,--\,\eqref{estvj} yield
\begin{eqnarray}\label{E-estvj}
 &&\hskip-4mm|v(t,x)|\le\sum\nolimits_{j}|v_j(t)|\le\sum\nolimits_j|v_j^0| e^{\,(\Phi-(\pi j/l)^2)\,t} \\
\nonumber
 && +\frac1{(\pi/l)^2-\Phi}
 \Big(\sum\nolimits_j e^{(1-\theta)(\Phi-(\pi j/l)^2)\,t}\,\frac{6}{\pi j}\Big)\sup\limits_{\,\tau\in[0,\,\theta t]}\nu(\tau)
 +\Big(\sum\nolimits_j\frac1{(\pi j/l)^2-\Phi}\cdot\frac{6}{\pi j}\Big)\sup\limits_{\,\tau\in[\theta t,\,t]}\nu(\tau).
\end{eqnarray}
 Hence, (\ref{dfftx}), (\ref{dfvtx}), \eqref{E-estvj} and Schwartz's inequality
imply the desired~\eqref{estdifference}.

\vskip1mm
$(ii)$
Since $\tilde\mu_0=\tilde\mu_1=0$, we can choose $\tilde\varphi(x)\equiv 0$ as
a~solution of the stationary problem (\ref{statprob}).
Hence, $v=\varphi(t,x)-U(t,x)$ and $v_0=\varphi_0(x)-U(0,x)$, see \eqref{dfvtx} and \eqref{dfv0U},
and \eqref{boundprv} takes the form
\begin{equation}\label{boundprv-2}
 \partial_tv=v_{\,,xx}+(\pi/l)^2 v+f(t,x),\quad v(0,x)=v_0(x),\quad v(t,0)=v(t,l)=0,
\end{equation}
where $f=(\pi/l)^2U-\dt U$.
As~in the proof of claim $(i)$, one may represent a~solution of \eqref{boundprv-2} in the series form
$v(t,x)=\sum\nolimits_{j=1}^\infty v_j(t)\sin(\pi j\,x/l)$,
 where $v_j(t)$ solves (\ref{Cauchvj}) with $\Phi=(\pi/l)^2$.

For $j=1$ we obtain the Cauchy's problem
\[
 v_1^{\,\prime}=f_1(t),\qquad v_1(0)=v_1^0,
\]
where $f_1(t)=\frac{2}{\pi}\,\big((\pi/l)^2(\delta_0(t)+\delta_1(t)) -\delta'_0(t)-\delta'_1(t)\big)$,
see (\ref{dfftx}), and
\[
  v_1^0 = \frac2\pi \int_0^l v_0(s)\sin(\pi s/l)\ds
 =\frac{2}{l}\int_0^l \varphi_0(s)\sin(\pi s/l)\ds
 -\frac2\pi\,\big(\delta_0(0)+\delta_1(0)\big).
\]
 Hence, $v_1(t)=v_1^0+\int_0^\infty\!f_1(\tau)\dtau-\int_t^\infty \!f_1(\tau)\dtau$,
where the improper integrals converge in view of~(\ref{addcond}).
By \eqref{E-fseries} for $j=1$, we obtain
\[
 \big|v_1(t)-v_1^0-\int_0^\infty f_1(\tau)\dtau\big|\le \frac{6}{\pi}\int_t^\infty\nu(\tau)\dtau.
\]
For $j>1$ we get the Cauchy's problem
\[
 v_j^{\,\prime}=(\pi/l)^2(1 - j^2) v_j +f_j(t),\qquad v_j(0)=v_j^0,
\]
 By Lemma~\ref{L-stabODE} with $a=(\pi/l)^2(1-j^2)$ and $\nu(t)=f_j(t)$,
 we get the estimates for $v_j\ (j\ge 2)$, see~(\ref{estvj}),
\begin{equation}\label{estvj-2}
 |v_j(t)|\le|v_j^0| e^{\,(\pi/l)^2(1-j^2)\,t} +(\pi/l)^{-2}(j^{\,2}-1)^{-1}\big(
 e^{(1-\theta)(\pi/l)^2(1-j^2)\,t}\!\!\sup\limits_{\,\tau\in[0,\,\theta t]}|f_j(\tau)|
 +\!\!\sup\limits_{\,\tau\in[\theta t,\,t]}|f_j(\tau)|\big)
\end{equation}
for any $\theta\in(0,1)$. Hence, see \eqref{E-fseries} and \eqref{E-estvj},
\begin{eqnarray}\label{E-estvj-2}
 &&\hskip-18mm
 \sum\nolimits_{j\ge2}|v_j(t)|
 \le \sum\nolimits_{j\ge2}|v_j^0|\,e^{\,(\pi/l)^2(1-j^2)\,t}
 +\Big(\sum\nolimits_{j\ge2}\frac{(l/\pi)^2}{j^2-1}\cdot\frac{6}{\pi j}\Big)
 \sup\nolimits_{\,\tau\in[\theta t,\,t]}\nu(\tau)\\
\nonumber
 &&+\Big(\sum\nolimits_{j\ge2}\frac{(l/\pi)^2}{j^2-1}\, e^{\,(1-\theta)(\pi/l)^2(1-j^2)\,t}
 \,\frac{6}{\pi j}\Big)\sup\nolimits_{\,\tau\in[0,\,\theta t]}\nu(\tau)\\
\nonumber
 &&
 \le\big(\sum\nolimits_{j\ge2} e^{\,2(\pi/l)^2(1-j^2)\,t}\big)^{1/2}\|v_0\|_{L_2}\\
 && +\,\frac{2\,l^2}{\pi^3}\sum\nolimits_{j\ge2} e^{\,(1-\theta)(\pi/l)^2(1-j^2)\,t}
 \,\sup\nolimits_{\,\tau\in[\theta t,\,t]}\nu(\tau)
\nonumber
  +\frac{3\,l^2}{2\,\pi^3}\sup\nolimits_{\,\tau\in[0,\,\theta t]}\nu(\tau).
\end{eqnarray}
Here we used $\sum\nolimits_{j\ge2}\frac1{j(j^2-1)}=\frac14$. Note that
\[
 \varphi(t,x)-\varphi_\infty(x)
 =\!\big(v_1(t)-v_1^0-\!\int_0^\infty\!\!f_1(\tau)\dtau\big)\sin(\pi x/l)
 +\!\sum\nolimits_{j\ge2} v_j(t)\sin(\pi j x/l)+U(t,x).
\]
Denote by
\begin{eqnarray*}
 \widetilde M_0(t)\eq\Big[\sum\nolimits_{j=2}^\infty e^{\,2\,(\pi/l)^2(4-j^2)\,t}\Big]^{1/2},\qquad
 \widetilde M_1(t)=\sum\nolimits_{j=2}^\infty e^{\,(1-\theta)\,(\pi /l)^2(4-j^2)\,t}\,.
\end{eqnarray*}
By the above, we have
\begin{eqnarray}\label{estdifference2}
 &&\hskip-7mm
 \big|\varphi(t,x)-\varphi_\infty(x)\big|
 \le \frac6\pi\int_t^\infty\nu(\tau)\dtau +\max\{|\delta_0(t)|,|\delta_1(t)|\}
 +\,\widetilde M_0(t) e^{\,-3(\pi/l)^2 t}\Vert v_0\Vert_{L_2}
 \nonumber\\
 &&\hskip23mm
 +\,\widetilde M_1(t)\frac{2\,l^2}{\pi^3}\,
 e^{\,-3\,(1-\theta)(\pi/l)^2 t}\sup\nolimits_{\,\tau\in[0,\,\theta t]}\nu(\tau)
 +\frac{3\,l^2}{2\,\pi^3}\,\sup\nolimits_{\,\tau\in[\theta t,\,t]}\nu(\tau)
\end{eqnarray}
for any $\theta\in(0,1)$ and all $(t,x)\in[0,\infty)\times[0,l]$.
This estimate implies the claim~$(ii)$.\qed

\begin{remark}\rm In Theorem~\ref{T-stabnonstat}, the limit metric $\tilde g$ has $r=\Phi\,\hat g$, which
may be as positive so negative definite.
In~case $(i)$, the system has a global single point attractor,
while in case $(ii)$ the limit solution metric depends on initial condition.
\end{remark}

\begin{remark}\rm
For $\Phi=0$ and $\mu_k(t)\equiv\tilde\mu_k$ the function $\tilde\varphi(x)=\tilde\mu_1(l-x)+\tilde\mu_0\,x$ is linear, and we have $\delta_k(t)=0,\,U=0,\,\nu(t)=0$, $v_0=\varphi_0-\tilde\varphi(x)$. Hence \eqref{estdifference} reads as
$|\varphi(t,x)-\tilde\varphi(x)|\le M_0\|\varphi(t,\cdot)-\tilde\varphi\|_{L_2}$.
\end{remark}

\begin{lemma}\label{L-stabODE}
Let $y(t)$ solves the Cauchy's problem $($for the ODE$)$
\[
 y^\prime=\alpha(t)y+\nu(t),\qquad
 y(0)=y_0,
\]
where
the functions $\alpha,\nu\in C[0,\infty)$,
$\alpha(t)\le a<0$ and $\nu(t)$ is bounded.
Then
\begin{equation}\label{estsolODE}
 |y(t)|\le|y_0| e^{\,a t}+|a|^{-1}e^{(1-\theta)\,a t}\sup\nolimits_{\,\tau\in[0,\,\theta t]}|\nu(\tau)|
 +|a|^{-1}\sup\nolimits_{\,\tau\in[\theta t,\,t]}|\nu(\tau)|
\end{equation}
for any $\theta\in(0,1)$. In particular, if $\,\lim_{\,t\rightarrow\infty}\nu(t)=0$ then
$\lim_{\,t\rightarrow\infty}y(t)=0$.
\end{lemma}

\noindent\textbf{Proof}.
As is known, $y(t)=y_0 e^{\,\int_0^t\alpha(t)\,{\rm dt}}
+\int_0^t e^{\,\int_\tau^t\alpha(\tau)\,{\rm d\tau}}\nu(\tau)\dtau$. Hence we have the estimate
\begin{eqnarray*}
 |y(t)|\eq|y_0| e^{\,at}+\int_0^{\theta t} e^{\,a(t-\tau)}|\nu(\tau)|\dtau
 +\int_{\theta t}^t e^{\,a(t-\tau)}|\nu(\tau)|\dtau\\
 &\le&|y_0| e^{\,at}+\sup\nolimits_{\,\tau\in[0,\,\theta t]}|\nu(\tau)|\int_0^{\theta t} e^{\,a(t-\tau)}\dtau
 +\sup\nolimits_{\,\tau\in[\theta t,\,t]}|\nu(\tau)|\int_{\theta t}^t e^{\,a(t-\tau)}\dtau.
\end{eqnarray*}
The above and
$\int_0^{\theta t} e^{\,a(t-\tau)}\dtau=(e^{\,at}-e^{(1-\theta)\,at})/a$,
\
$\int_{\theta t}^t e^{\,a(t-\tau)}\dtau=(e^{(1-\theta)\,at}-1)/a$ yield~(\ref{estsolODE}).
\qed

\begin{example}\rm
Consider the problem (\ref{E-r}) with $\Phi=(\pi/l)^2$ and stationary boundary conditions,
$\mu_j(t)\equiv\tilde\mu_j\ (j=0,1)$. Assume that
$\tilde\mu_0\neq -\tilde\mu_1$, then the stationary boundary problem (\ref{statprob})
does not have a solution (see Lemma~\ref{Ex-tilder}). Let us show that in this case the solution of
the non-stationary problem (\ref{E-r}) is not stabilized for
$t\rightarrow\infty$. As above, if $\varphi(t,x)$ is the solution of (\ref{E-r}),~then
 $v(t,x)=\varphi(t,x)-U_0(x)$ with $U_0(x)=\frac{l-x}{l}\,\tilde\mu_0+\frac{x}{l}\,\tilde\mu_1$
is the solution of the problem
\[
 \partial_tv=v_{\,,xx}+(\pi/l)^2v+f(x),\quad v(0,x)=v_0(x),\quad v(t,0)=v(t,l)=0,
\]
where $f(x)=(\pi/l)^2U_0(x)$ and $v_0(x)=\varphi_0(x)-U_0(x)$.
Applying the method used in the proof of
Theorem~\ref{T-stabnonstat}, we obtain after easy calculations:
\begin{eqnarray*}
 v(t,x)\eq(v_1^0+f_1t)\sin(\pi x/l)\\
 \plus\sum\nolimits_{j=2}^\infty\Big(\big(v_j^0-\frac{f_j}{(\pi/l)^2(j^2-1)}\big)
 e^{\,(\pi/l)^2(1-j^2)\,t}
 +\frac{f_j}{(\pi/l)^2(j^2-1)}\Big)\sin(\pi j\,x/l),
\end{eqnarray*}
where $v_j^0$ and $f_j$ are Fourier's coefficients of the functions
$v_0(x)$ and $f(x)$, respectively. In particular,
$f_1=\frac{2n\pi^2}{l^3}(\tilde\mu_0+\tilde\mu_1)\int_0^ls\,\sin(\pi s/l)\ds\neq
0$ (since $\tilde\mu_0+\tilde\mu_1\neq 0$).
Hence $|v(t,x)|\rightarrow\infty$ as $t\rightarrow\infty$, that is
$v(t,x)$ is not stabilized at infinity.
\end{example}

\end{document}